\pdfoutput=1
\documentclass[11pt]{article}
\usepackage[utf8]{inputenc}
\usepackage[T1]{fontenc}
\usepackage[british]{babel}
\usepackage{amsmath,amsfonts,amssymb,amsthm}
\usepackage{mathtools}
\usepackage{booktabs}
\usepackage{graphicx}
\usepackage{xcolor}
\usepackage[round,authoryear]{natbib}
\usepackage[a4paper,margin=2.6cm]{geometry}
\usepackage{hyperref}
\hypersetup{colorlinks=true, linkcolor=blue, citecolor=blue, urlcolor=blue}
\allowdisplaybreaks
\usepackage{authblk}

\newtheorem{proposition}{Proposition}
\theoremstyle{remark}

\newcommand{\Bb}{\bar{B}}
\newcommand{\din}{\delta^{\mathrm{in}}}
\newcommand{\dout}{\delta^{\mathrm{out}}}
\newcommand{\Jbar}{\bar{J}}
\newcommand{\LB}{\mathrm{LB}}
\newcommand{\UB}{\mathrm{UB}}

\usepackage{setspace}
\title{Ride-time inequalities for the event-based formulation of the dial-a-ride problem}
\author{Julian Golak}
\author{Arne Schulz}

\affil{\footnotesize Institute of Operations Management, University of Hamburg Business School, Hamburg, Germany.\\ \texttt{\{julian.golak, arne.schulz\}@uni-hamburg.de}}

\date{\today}

\begin{document}
\maketitle

\begin{abstract}
The paper investigates the model structure of the event-based mixed-integer programming formulation of the Dial-a-Ride Problem. We find that the design of events allows introducing new valid inequalities on the ride times. Time progression is modelled using Big $M$ constraints which lead to a weak LP relaxation and deteriorate the solution process. The newly introduced valid inequalities deal with this issue and allow reducing the computation times by up to 99\% for benchmark instances (a-X$^\star$ instances) in comparison to the state-of-the-art event-based model formulation. To evaluate the potential of the new valid inequalities, we increase the time windows on the benchmark set by another 15 minutes and evaluate the model on the same instances with a time-related objective function where the Big-$M$-influenced time variables directly occur in the objective function. Even for the large time windows only two instances could not be solved to optimality within two hours (final gaps of at most 0.8\%). With the time-oriented objective the largest average computation time with our valid inequalities was 74 seconds.
\end{abstract}


\section{Introduction} \label{sec:introduction}

The Dial-a-Ride Problem (DARP) searches for optimal tours transporting passengers from their pickup to their delivery location while respecting vehicle capacities, time windows, and maximum ride time constraints. The traditional application is the transportation of elderly or disabled persons, but nowadays it has also applications in ridepooling services for the general public \citep{pfeiffer2022alns}. The DARP was introduced by \citet{psaraftis1980dynamic}.

With respect to exact optimisation procedures \citet{psaraftis1980dynamic} presented a dynamic program for the case with one vehicle. Later, the focus shifted to mixed-integer programming formulations starting with \citet{cordeau2006branch} who presented a 3-index formulation. The model was advanced to a 2-index formulation by \citet{ropke2007models} for homogeneous vehicles. Both model formulations use binary variables picturing whether a vehicle moves directly between two locations $i$ and $j$. \citet{gaul2025tight} call these models location-based formulations.

More recently, \citet{rist2021new} developed a new type of model formulation based on so-called fragments. A fragment is a part of the tour where the vehicle is empty at the beginning and at the end. A vehicle's tour is then composed of fragments.

Only shortly after, \citet{gaul2022event} presented a further model formulation based on events. An event pictures the vehicles' current occupancy as a $Q$-tuple with the vehicle capacity $Q$. The detailed construction is presented in Section \ref{sec:graph}. Later on, \citet{gaul2025tight} improved their model and added new valid inequalities using the structure of  the event-based graph to improve the search within a branch-and-cut procedure.

Since then, event-based models have been applied to other problem settings in routing and beyond. \citet{stallhofer2025event} and \citet{zhao2026new} presented event-based models for the electric autonomous DARP. Further, it was used by \citet{anzenhofer2025potential} to evaluate the potential of governmental regulations on shared mobility-on-demand systems and by \citet{zhao2026dial} for the DARP with synchronized visits. Moreover, the event-based DARP model formulation was applied to the line-based DARP with \citep{barth2025line} and without transfers \citep{reiter2024line,lauerbach2026solving}.

\citet{he2025dynamic} suggest to apply the event-based structure for electric vehicle fleet management for ride-hailing; \citet{ferreira2025new} do the same for the Prisoner Transportation Problem. \citet{chi2026event} developed an event-based model for an inland multi-size container transportation problem.

These applications do not only show the potential of event-based DARP model formulations, but also highlight the need for further improvements. While the approach is naturally limited to small vehicle sizes, which however are typical in practice, the computational study by \citet{gaul2025tight} also showed that increasing time windows are a drawback. This is true for two reasons: First, larger time windows lead to an increased event graph. Second, larger time windows lead to larger Big $M$s in the constraints modelling the time progression while Big $M$ terms deteriorate the solution process \citep{schulz2024using}.

The paper at hand presents new valid inequalities to address the second issue. Therefore, we introduce the DARP formally in Section \ref{sec:problem}. Afterwards, the event graph is introduced in Section \ref{sec:graph}. Next, the event-based model formulation is presented in Section \ref{sec:formulations} including the newly introduced ride-time inequalities (Sections \ref{sec:ridecut} and \ref{sec:strengthened}). The resulting model formulation is evaluated on benchmark instances in a comprehensive computational study in Section \ref{sec:study}. The paper closes with a conclusion in Section \ref{sec:conclusion}.
\section{Problem formulation}
\label{sec:problem}

We consider the Dial-a-Ride Problem of \citet{cordeau2006branch} in the notation of \citet{gaul2025tight}.
A set of requests $R = \{1, \dots, n\}$ is served by at most $K$ identical vehicles of capacity $Q$ that start and end their tours at a common depot~$0$.
Request $i$ asks for the transport of $q_i$ passengers from its pickup location $i^+ = i$ to its delivery location $i^- = n + i$.
We write $P = \{1, \dots, n\}$ for the pickup locations, $D = \{n + 1, \dots, 2n\}$ for the delivery locations, $J = P \cup D$ and $\Jbar = J \cup \{0\}$.
Every location $j \in \Jbar$ has a service time $s_j$ and a time window $[e_j, \ell_j]$ for the start of the service; the depot has $s_0 = 0$ and $[e_0, \ell_0] = [0, T]$, where $T$ is the planning horizon.
Driving from location $i$ to location $j$ takes the time $t_{ij}$ and costs $c_{ij}$.
A vehicle may wait at a location before it starts the service.
The ride time of request $i$ is the time between the end of the service at $i^+$ and the start of the service at $i^-$; it is limited by $L$.
In the benchmark instances, only one side of a request carries a given time window.
Following \citet{cordeau2006branch}, the window of the other side is derived from it: for an outbound request with a given delivery window we set $e_{i^+} = \max\{0, e_{i^-} - L - s_{i^+}\}$ and $\ell_{i^+} = \ell_{i^-} - t_{i^+ i^-} - s_{i^+}$, and for an inbound request with a given pickup window we set $e_{i^-} = e_{i^+} + s_{i^+} + t_{i^+ i^-}$ and $\ell_{i^-} = \ell_{i^+} + s_{i^+} + L$.
Thus, $e_{i^-}$ is the earliest possible start of the delivery of request $i$ in every instance.

\section{Event graph}
\label{sec:graph}

The event-based model formulation of \citet{gaul2022event} describes the state of a vehicle after the service at a location by an event.
An event is a tuple $v = (v_1, \dots, v_Q) \in V$, where $v_1 \in J$ is the location that has just been served and $v_2 > \dots > v_Q$ are the requests on board besides the request of $v_1$, padded with zeros.
The request of $v_1$ is on board after $v$ if and only if $v_1$ is a pickup location.
The set of events $V$ consists of the depot event $v_0 = (0, \dots, 0)$ and of every tuple whose load respects the capacity $Q$ and whose requests are pairwise compatible, that is, every two requests on board admit a sequence of their four locations in the required order that respects the time windows, the ride-time limit and the horizon.
For a location $j \in J$ let $V_j$ be the events with $v_1 = j$, so that $V_{i^+}$ and $V_{i^-}$ are the events in which request $i$ is picked up and delivered.

The arc set $A$ connects the depot event with every event $(i^+, 0, \dots, 0)$ and every event $(i^-, 0, \dots, 0)$ with the depot event.
It connects an event $v$ with an event $w$ if $w_1$ is the pickup location of a request that is not on board after $v$ and for which a seat is free, or the delivery location of a request that is on board after $v$, and if the requests on board after $w$ are those on board after $v$ with that request added or removed.
Only arcs whose requests are pairwise compatible in the resulting order are kept, as described in \citet[Section~2.2]{gaul2025tight}.
An arc $a = (v, w)$ has the cost $c_a = c_{v_1 w_1}$ and the travel time $t_a = t_{v_1 w_1}$, and $\din(v)$ and $\dout(v)$ denote the arcs entering and leaving $v$, respectively.
By construction, every path from the depot event back to the depot event is a vehicle tour that respects the capacity and visits the pickup event of each request before its delivery event.
The time windows, the ride-time limit and the horizon are not part of the graph and are enforced by the formulation below.

\paragraph{Preprocessing.}
\citet[Section~4.1]{gaul2025tight} compute for every event $v$ a lower bound $B^{\LB}_v$ and an upper bound $B^{\UB}_v$ on the start of the service at $v_1$ in event $v$.
The lower bound propagates the earliest start of the service along the predecessors of $v$, and the upper bound propagates the latest start along the successors, including the ride-time limit of the requests on board.
We compute the bounds exactly as in their recursions, as a fixed point over the whole graph.
An event with $B^{\UB}_v < B^{\LB}_v$ and an arc $(v, w)$ with $B^{\LB}_v + s_{v_1} + t_{(v, w)} > B^{\UB}_w$ cannot be part of a feasible solution and are removed from the graph.
We refer to this graph preprocessing as GP.

\section{Event-based model formulation}
\label{sec:formulations}

\subsection{Location-augmented event-based model formulation}
\label{sec:laeb}

The location-augmented event-based model formulation (LAEB) of \citet{gaul2025tight} uses a binary variable $x_a$ for every arc $a \in A$ and one continuous variable $\Bb_j \ge 0$ for every location $j \in \Jbar$, the start of the service at $j$.
The variable $\Bb_0$ bounds the return times of all vehicles.
For two locations $i, j \in J$ let $A_{ij} = \{(v, w) \in A : v_1 = i, w_1 = j\}$ be the arcs from $i$ to $j$, and let $M_{ij} = \ell_i + s_i + t_{ij} - e_j$.
The formulation reads
\begin{align}
  \min\ & \sum_{a \in A} c_a x_a \label{eq:laeb:obj} \\
  \text{s.t.}\ & \sum_{a \in \din(v)} x_a - \sum_{a \in \dout(v)} x_a = 0 && \forall v \in V, \label{eq:laeb:flow} \\
  & \sum_{v \in V_{i^+}} \sum_{a \in \din(v)} x_a = 1 && \forall i \in R, \label{eq:laeb:serve} \\
  & \sum_{a \in \dout(v_0)} x_a \le K, \label{eq:laeb:fleet} \\
  & \Bb_j \ge \Bb_i + s_i + t_{ij} - M_{ij} \Bigl(1 - \sum_{a \in A_{ij}} x_a\Bigr) && \forall i, j \in J,\ A_{ij} \neq \emptyset, \label{eq:laeb:time} \\
  & \Bb_j \ge e_0 + t_{0j}\, x_{(v_0, (j, 0, \dots, 0))} && \forall j \in P, \label{eq:laeb:start} \\
  & \Bb_0 \ge \Bb_i + s_i + t_{i0} - M_{i0} \bigl(1 - x_{((i, 0, \dots, 0), v_0)}\bigr) && \forall i \in D, \label{eq:laeb:return} \\
  & e_j \le \Bb_j \le \ell_j && \forall j \in \Jbar, \label{eq:laeb:tw} \\
  & \Bb_{i^-} - \Bb_{i^+} - s_{i^+} \le L && \forall i \in R, \label{eq:laeb:ride} \\
  & x_a \in \{0, 1\} && \forall a \in A. \label{eq:laeb:bin}
\end{align}
The objective~\eqref{eq:laeb:obj} minimises the routing cost.
Constraints~\eqref{eq:laeb:flow} are the flow conservation constraints at every event, Constraints~\eqref{eq:laeb:serve} pick up every request exactly once, and Constraint~\eqref{eq:laeb:fleet} limits the number of tours.
Constraints~\eqref{eq:laeb:time} link the service start times of two locations whenever a vehicle drives from one to the other, with the usual Big $M$ relaxation otherwise; Constraints~\eqref{eq:laeb:start} and~\eqref{eq:laeb:return} do the same for the first and the last location of a tour.
Constraints~\eqref{eq:laeb:tw} are the time windows and Constraints~\eqref{eq:laeb:ride} the ride-time limit.
With the preprocessing bounds, \citet{gaul2025tight} replace the time windows~\eqref{eq:laeb:tw} for $j \in J$ by the tighter bounds
\begin{align}
  \sum_{v \in V_j} B^{\LB}_v \sum_{a \in \din(v)} x_a \le \Bb_j \le \sum_{v \in V_j} B^{\UB}_v \sum_{a \in \din(v)} x_a && \forall j \in J, \label{eq:laeb:gp}
\end{align}
which select the bounds of the event through which location $j$ is visited; we keep~\eqref{eq:laeb:tw} in addition.

\subsection{Valid inequalities from literature}
\label{sec:gaulineq}

\citet[Section~4.2]{gaul2025tight} propose three families of valid inequalities for the event-based models: the infeasible-path inequalities, which exclude pairs of consecutive arcs and short paths that no feasible schedule can use given the preprocessing bounds; the vehicle-sharing inequalities, which state that two requests either share a vehicle or are served one after the other; and the customer-incompatibility inequalities for pairs of requests that cannot share a vehicle.
We use all families with the dominance rules of their implementation, in the combination of their final benchmark, and refer to the formulation with GP and all families as the base configuration.
The branch-and-cut node methods of their Sections~4.3--4.5 are not used, as in their final benchmark.

\subsection{Ride-time inequalities}
\label{sec:ridecut}

The time propagation constraints~\eqref{eq:laeb:time} are the main link between the arc variables and the time variables, and their Big $M$ relaxation becomes weak when $\Bb_i < \ell_i$ or $\Bb_j > e_j$. Then, the constraint is satisfied by a positive $\sum_{a\in A_{ij}} x_a$ although $\Bb_j < \Bb_i + s_i + t_{ij}$, i.e. the arcs from $i$ to $j$ can carry flow in the LP relaxation without the time variables reflecting it. The reason is the Big $M$ and its choice, which is, however, necessary to ensure feasibility in the extreme case with $\Bb_i = \ell_i$ and $\Bb_j = e_j$.

In this section, we use the structure of the event graph to overcome the drawback of the Big $M$ formulation. The approach follows the idea of fixed paths from \citet{schulz2024branch,schulz2024using}. Here we do not need a completely fixed path, but we know a fixed sequence of events if they are used.

Due to construction of the event graph, we know that location $i^-$ is visited after location $i^+$ in the same tour. Thus, we do not need a Big $M$ for this relation, as $\Bb_{i^-} > \Bb_{i^+}$ certainly holds. Moreover, we know that all events with request $i$ on positions $2, ..., Q$ in the event must either be visited between location $i^+$ and location $i^-$ or are not used in the solution.

Thus, we get a path
\begin{equation} \label{eq:path}
(i^+, ...) \rightarrow (b_1, ..., i, ...) \rightarrow (b_2, ..., i, ...) \rightarrow ... \rightarrow (b_k, ..., i, ...) \rightarrow (i^-, ...)
\end{equation}
between the pickup and delivery location of request $i$. If this path is chosen, it follows
\begin{equation} \label{eq:path1}
\Bb_{i^-} \geq (\Bb_{i^+} + s_{i^+}) + (t_{i^+ b_1} + s_{b_1}) + (t_{b_1 b_2} + s_{b_2}) + \dots + t_{b_k i^-} = \Bb_{i^+} + \sum_{l=0}^{k} \bigl(s_{b_l} + t_{b_l b_{l+1}}\bigr)
\end{equation}
with $b_0 = i^+$ and $b_{k+1} = i^-$. This path is chosen if the corresponding arc variables are 1.

Generally, we formulate inequalities linking the ride time of a request directly to the arcs on which it is on board, without a Big $M$ term as follows. For a request $i \in R$ let
\begin{align*}
  A_i = \bigl\{(v, w) \in A : v_1 = i^+ \text{ or } i \in \{v_2, \dots, v_Q\}\bigr\}
\end{align*}
be the arcs on which request $i$ is on board, that is, the arcs leaving an event after which $i$ sits in the vehicle.
In every solution, the arcs of $A_i$ with $x_a = 1$ form exactly the path from the pickup event of $i$ to the delivery event of $i$.
Summing Constraints~\eqref{eq:laeb:time} along this path gives
\begin{align}
  \Bb_{i^-} - \Bb_{i^+} \ge \sum_{a = (v, w) \in A_i} \bigl(s_{v_1} + t_a\bigr) x_a && \forall i \in R. \label{eq:ridecut}
\end{align}

\begin{proposition} \label{prop1} \ \\
Inequalities \eqref{eq:ridecut} are valid for model \eqref{eq:laeb:flow}--\eqref{eq:laeb:bin}.
\end{proposition}

\begin{proof}
In an integer feasible solution, the graph construction ensures with Constraints \eqref{eq:laeb:flow}, \eqref{eq:laeb:serve} and \eqref{eq:laeb:time} (no cycles) that a single path of the structure of \eqref{eq:path} exists. For it, \eqref{eq:path1} and by it \eqref{eq:ridecut} holds.
\end{proof}

\begin{proposition} \label{prop2} \ \\
Valid inequalities \eqref{eq:ridecut} ensure for model \eqref{eq:laeb:flow}--\eqref{eq:laeb:bin} that ride time feasibility holds in the weighted flow average in the LP solution.
\end{proposition}

\begin{proof}
Due to the graph construction and flow conservation constraints \eqref{eq:laeb:flow} every flow leaving a pickup event of request $i$ (i.e. with $v_1 = i^+$) transfers to a delivery event of request $i$ (i.e. with $w_1 = i^-$). Moreover, a total flow of 1 has to leave the pickup events of every request $i$ due to \eqref{eq:laeb:serve}. In a fractional solution, the entire flow of size 1 can be decomposed into several flows with a size between 0 and 1. For each of these fractional flows between pickup and delivery location of a request $i$ a path \eqref{eq:path} exists such that \eqref{eq:ridecut} sums up the service and travel times on the right side, i.e. aggregates the fractional flows including possible cycles, which have nonnegative length, within these paths.

Together with the ride-time limit~\eqref{eq:laeb:ride}, valid inequalities \eqref{eq:ridecut} imply
\[
\sum_{a \in A_i} (s_{v_1} + t_a) x_a \le L + s_{i^+},
\]
a constraint on the arc variables alone that cuts off fractional solutions in which the arcs carrying request $i$ are in the weighted flow average longer than the ride-time limit allows.
\end{proof}

Note that without valid inequalities \eqref{eq:ridecut} ride time feasibility in the sense of Proposition~\ref{prop2} is not guaranteed in the LP relaxation due to the Big-$M$-type Constraints \eqref{eq:laeb:time}.

In the pure event-based formulation of \citet{gaul2022event}, where the time variables $B_v$ are attached to events $v \in V$, the inequality has a direct counterpart, when $B$ variables of not used events are set to 0
\begin{align} \nonumber
0 \leq B_v \leq B_v^{UB} \sum_{a \in \delta^{in}(v)} x_a && \forall v \in V.
\end{align}
Then, $\Bb_{i^-} - \Bb_{i^+}$ can be replaced in \eqref{eq:ridecut} by
\[
\sum_{v \in V_{i^-}} B_v - \sum_{v \in V_{i^+}} B_v.
\]
However, this type of constraint cannot be transferred to further relations between pickup and delivery locations, as we do not know whether different requests are served by the same vehicle.

\subsection{Strengthened ride-time inequalities} \label{sec:strengthened}
Valid inequalities \eqref{eq:ridecut} can be strengthened by replacing $s_{v_1} + t_a$ by
\[
\tau_a = \max\{ s_{v_1} + t_a, B^{LB}_w - B^{UB}_v \}.
\]
This leads to inequalities
\begin{align} \label{eq:lridecut}  \Bb_{i^-} - \Bb_{i^+} \ge \sum_{a = (v, w) \in A_i} \tau_a \cdot x_a && \forall i \in R.
\end{align}

\begin{proposition} \ \\
Inequalities \eqref{eq:lridecut} are valid for model \eqref{eq:laeb:flow}--\eqref{eq:laeb:bin}.
\end{proposition}

\begin{proof}
As Proposition \ref{prop1} shows that \eqref{eq:ridecut} is valid, it remains to show that replacing $s_{v_1} + t_a$ by $B^{LB}_w - B^{UB}_v$ does not violate feasibility. Thus, we assume that
\[
B^{LB}_w - B^{UB}_v > s_{v_1} + t_a.
\]
If the vehicle starts service at location $v_1$ at the latest possible time $B^{UB}_v$ and drives directly to location $w_1$, the vehicle still arrives before the tightened time window beginning $B^{LB}_w$. If service at $v$ starts before $B^{UB}_v$, waiting time at $w$ further increases. Hence, service cannot start by arrival but after a waiting time of at least
\[
B^{LB}_w - B^{UB}_v - s_{v_1} - t_a
\]
which is exactly the amount by which $\tau_a$ exceeds $s_{v_1} + t_a$ on the right-hand side of the inequality in comparison to \eqref{eq:ridecut}. Thus, \eqref{eq:lridecut} is valid.
\end{proof}

The proof shows that \eqref{eq:lridecut} strengthens \eqref{eq:ridecut} by necessary waiting time. Proposition \ref{prop4} holds analogously to Proposition \ref{prop2}.

\begin{proposition} \label{prop4} \ \\
Valid inequalities \eqref{eq:lridecut} ensure for model \eqref{eq:laeb:flow}--\eqref{eq:laeb:bin} that ride time feasibility holds in the weighted flow average in the LP solution.
\end{proposition}

\begin{proof}
Analogous to Proposition \ref{prop2}.
\end{proof}

Note that valid inequalities \eqref{eq:lridecut} dominate valid inequalities \eqref{eq:ridecut} since $\tau_a \geq s_{v_1} + t_a$.

\subsection{Objectives}
\label{sec:objectives}
Besides the routing cost~\eqref{eq:laeb:obj}, we study a customer-oriented objective in which the time variables appear directly in the objective, so that the time constraints have to be tight in every solution.
We use the relative detour of \citet{pfeiffer2022alns,pfeiffer2022new}, the delay of every passenger against the earliest possible delivery, measured in units of the direct ride time:
\begin{align}
  \min\ \sum_{i \in R} q_i\, \frac{\Bb_{i^-} - e_{i^-}}{t_{i^+ i^-}}. \label{eq:detour}
\end{align}
A request that is delivered at its earliest possible time contributes zero.
\citet{pfeiffer2022alns,pfeiffer2022new} measure the detour against $e_{i^+} + t_{i^+ i^-}$, which equals $e_{i^-}$ in their setting without service times.
On the benchmark instances the two readings have the same optimal solutions and differ by the instance constant
\begin{align*}
  \sum_{i \in R} q_i\, \frac{e_{i^-} - e_{i^+} - t_{i^+ i^-}}{t_{i^+ i^-}},
\end{align*}
which is the service time at the pickup for inbound requests and the slack between the earliest pickup and the delivery window for outbound requests, each in units of the direct ride time.
Objective~\eqref{eq:detour} is used with LAEB, where the delivery times are location variables; all constraints remain as in Section~\ref{sec:laeb}. Results with this objective function are presented in Appendix \ref{sec:results:detour}.
\section{Computational study}
\label{sec:study}
\paragraph{Setting.}
All models are implemented in Python~3.12 with Gurobi~13.0.2 and solved on the Hummel-2 cluster of the University of Hamburg on nodes with two AMD EPYC 9654 processors (96 cores each, 2.4~GHz) and 768~GB of RAM, with one thread per solve, a relative optimality gap of zero and a time limit of 7200 seconds, as in \citet{gaul2025tight}.
The CPU times in the tables are the pure Gurobi times.
Building the event graph with the preprocessing and the model takes 0.4 to 3.6 seconds per instance on average, depending on the set, and at most 23 seconds (a8-96-X$^{\star\star}$), and is the same for the three configurations up to the ride-time inequalities.
Every reported solution is verified independently of the solver: the tours are rebuilt from the arc variables, and the routing cost, the time windows, the ride times, the capacity, the precedence of pick-up and delivery, the horizon and the fleet size are rechecked from the instance data.
The optimal values of all instances solved by both studies coincide with those of \citet{gaul2025tight}.

\paragraph{Configurations.}
We compare, for LAEB, the plain model~\eqref{eq:laeb:obj}--\eqref{eq:laeb:bin}, the model with graph preprocessing and the bounds~\eqref{eq:laeb:gp} (GP), and the ``base'' configuration with GP and all valid inequalities of Section~\ref{sec:gaulineq}, each without ride-time inequalities, with the ride-time inequalities~\eqref{eq:ridecut} and with the strengthened ride-time inequalities~\eqref{eq:lridecut}.
In the detailed tables, we compare base, base + \eqref{eq:ridecut} and base + \eqref{eq:lridecut}.

The base configuration and the two configurations with ride-time inequalities are run with three Gurobi seeds and reported as means. The Gurobi seed steers decisions (e.g. primal heuristics, tie breaking, cut selection) within the Gurobi Branch-and-Cut procedure which have an influence on the results with roughly a factor of 3. For example, the three runs for instance a7-84-X$^{\star\star}$ with base + \eqref{eq:ridecut} required 6268/2944/6935\,s and those for b8-96-X$^{\star\star}$ with base 898/261/213\,s. A single run would not be able to isolate the effect of the inequalities from this. We, therefore, report mean results over three independent runs in the tables.

Results for plain and GP are reported for the sake of completeness shortly in the text. The benefit of ride time inequalities is so strong here that we restricted ourselves to a single Gurobi seed/run.

\subsection{Test instances}
\label{sec:instances}

We use the four instance sets of the computational study of \citet{gaul2025tight}.
The sets a and b are the benchmark instances introduced by \citet{cordeau2006branch} and \citet{ropke2007models}; each set contains 24 instances with $K \in \{2, \dots, 8\}$ vehicles and $n \in \{16, \dots, 96\}$ requests, named $aK$-$n$ and $bK$-$n$.
In set a, every request has one passenger, every service takes 3 minutes, the capacity is $Q = 3$ and the ride-time limit is $L = 30$ minutes.
In set b, a request has $q_i \in \{1, \dots, 6\}$ passengers, the service time equals the number of passengers, $Q = 6$ and $L = 45$ minutes.
Half of the requests are outbound with a given delivery window and half are inbound with a given pickup window; every given window is 15 minutes long, and the other side is derived as described in Section~\ref{sec:problem}.
The planning horizon is $T = 60\, n / K$ minutes, that is, between 4 and 12 hours.
Locations are points in the plane; travel times and routing costs are the Euclidean distances rounded to two decimals.

The sets a-X$^\star$ and b-X$^\star$ are variants where the given time window of every request, that is, the delivery window of an outbound and the pickup window of an inbound request, is extended by 15 minutes at its end, so that it is 30 minutes long.
The window of the other side is derived from it and is thereby extended by the same amount.
All other data, including the capacity and the ride-time limit, are unchanged (adapted from \citet{gschwind2015effective} as in \citet{gaul2025tight}).
These sets are harder than the original ones, because wider windows admit more feasible combinations of requests per vehicle. Since every a-X$^\star$ and b-X$^\star$ instance is solved with the ride-time inequalities, we add the sets a-X$^{\star\star}$ and b-X$^{\star\star}$, in which every given window is extended by 30 minutes instead of 15, so that it is 45 minutes long; all other data are unchanged, and no reference values exist for them.

\subsection{Forced waiting}
\label{sec:results:wait}
Inequalities~\eqref{eq:lridecut} differ from~\eqref{eq:ridecut} only on arcs $a = (v, w)$ with $\tau_a > s_{v_1} + t_a$, on which a vehicle that starts the service at $v_1$ at the latest possible time $B^{\UB}_v$ still arrives at $w_1$ before $B^{\LB}_w$.
Table~\ref{tab:wait_stats} counts these arcs among the arcs with a passenger on board, $\bigcup_{i \in R} A_i$, on the preprocessed graphs.
Since the upper bounds of GP already propagate the time windows of the successors, forced waiting remains only where $\ell_{v_1}$ is tight and $e_{w_1}$ is late, and it is rare: 4.3 to 8.6~percent of the arcs with a passenger on board carry a positive waiting time, of 5.7 to 14.0 minutes on average and up to 67 minutes, and 53 to 83~percent of the requests have at least one such arc in $A_i$.
With the location windows $[e_{v_1}, \ell_{v_1}]$ in place of the preprocessing bounds, as in the plain model, the shares are 0.4 to 1.8 percentage points lower.
We therefore expect~\eqref{eq:lridecut} to differ little from~\eqref{eq:ridecut} in the computations although it is by definition stronger.

\begin{table}[htbp]
\centering\footnotesize
\caption{Arcs with forced waiting after the graph preprocessing, totals over the 24 instances of each set.}
\label{tab:wait_stats}
\begin{tabular}{l rrr rr r}
\toprule
Set & Arcs on board & With waiting & Share & Mean (min) & Max (min) & Requests \\
\midrule
a & 28192 & 1774 & 6.3\% & 5.7 & 24.7 & 684/1108 (62\%) \\
b & 26318 & 2256 & 8.6\% & 10.2 & 37.9 & 584/1108 (53\%) \\
a-X$^\star$ & 93991 & 4241 & 4.5\% & 8.1 & 38.5 & 806/1108 (73\%) \\
b-X$^\star$ & 63989 & 4198 & 6.6\% & 11.5 & 51.3 & 613/1108 (55\%) \\
a-X$^{\star\star}$ & 213749 & 9096 & 4.3\% & 10.9 & 54.4 & 923/1108 (83\%) \\
b-X$^{\star\star}$ & 121927 & 6394 & 5.2\% & 14.0 & 66.6 & 653/1108 (59\%) \\
\bottomrule
\end{tabular}
\par\vspace{3pt}
\begin{minipage}{\textwidth}\footnotesize
Notes: Arcs on board: arcs $(v, w)$ with at least one request on board after $v$, that is, $\bigcup_{i \in R} A_i$. With waiting: those with $\tau_a > s_{v_1} + t_a$; Share: their percentage. Mean, Max: mean and largest $\tau_a - s_{v_1} - t_a$ over these arcs in minutes. Requests: requests whose inequality~\eqref{eq:lridecut} differs from~\eqref{eq:ridecut}, over all requests of the set. Bounds $B^{\LB}$, $B^{\UB}$ of the graph preprocessing on the preprocessed graph.
\end{minipage}
\end{table}

\subsection{Routing cost}
\label{sec:results:cost}

Tables~\ref{tab:ridecut_LAEB_a_study}--\ref{tab:ridecut_LAEB_bX_study} compare the base configuration with the base configuration plus the ride-time inequalities~\eqref{eq:ridecut} and plus the strengthened inequalities~\eqref{eq:lridecut} on the four instance sets from literature.
For every instance they report the objective value and, for the three configurations, the gaps of the LP relaxation and of the root node, the number of branch-and-bound nodes and the CPU time as means over the three seeds, and the relative change of the CPU time against the base configuration; the notes below the tables define the columns.

\begin{table}[htbp]
\centering\scriptsize\setlength{\tabcolsep}{1.2pt}\renewcommand{\arraystretch}{0.9}
\caption{Base versus base + cuts on the a instances (routing costs~\eqref{eq:laeb:obj}).}
\label{tab:ridecut_LAEB_a_study}
\begin{tabular}{l r rrrr rrrrr rrrrr}
\toprule
 & & \multicolumn{4}{c}{base} & \multicolumn{5}{c}{base + \eqref{eq:ridecut}} & \multicolumn{5}{c}{base + \eqref{eq:lridecut}} \\
\cmidrule(lr){3-6}\cmidrule(lr){7-11}\cmidrule(lr){12-16}
Inst. & Obj. & LPgap & rGap & Nodes & CPU & LPgap & rGap & Nodes & CPU & $\Delta$CPU & LPgap & rGap & Nodes & CPU & $\Delta$CPU \\
\midrule
a2-16 & 294.3 & 0.4\% & 0.0\% & 1 & 0.03 & 0.4\% & 0.0\% & 1 & 0.02 & -- & 0.4\% & 0.0\% & 1 & 0.02 & -- \\
a2-20 & 344.9 & 2.1\% & 0.0\% & 1 & 0.04 & 1.3\% & 0.0\% & 1 & 0.07 & -- & 1.3\% & 0.0\% & 1 & 0.02 & -- \\
a2-24 & 431.1 & 2.3\% & 0.0\% & 1 & 0.05 & 2.1\% & 0.0\% & 1 & 0.05 & -- & 2.1\% & 0.0\% & 1 & 0.05 & -- \\
a3-18 & 300.5 & 0.8\% & 0.0\% & 1 & 0.03 & 0.8\% & 0.0\% & 1 & 0.04 & -- & 0.8\% & 0.0\% & 1 & 0.03 & -- \\
a3-24 & 344.8 & 2.4\% & 0.0\% & 1 & 0.09 & 2.2\% & 0.0\% & 1 & 0.08 & -- & 2.2\% & 0.0\% & 1 & 0.08 & -- \\
a3-30 & 494.8 & 1.8\% & 0.0\% & 1 & 0.04 & 1.4\% & 0.0\% & 1 & 0.05 & -- & 1.4\% & 0.0\% & 1 & 0.04 & -- \\
a3-36 & 583.2 & 2.5\% & 0.0\% & 1 & 0.19 & 2.5\% & 0.0\% & 1 & 0.08 & -- & 2.5\% & 0.0\% & 1 & 0.09 & -- \\
a4-16 & 282.7 & 0.7\% & 0.0\% & 1 & 0.12 & 0.6\% & 0.0\% & 1 & 0.07 & -- & 0.6\% & 0.0\% & 1 & 0.05 & -- \\
a4-24 & 375.0 & 1.7\% & 0.0\% & 1 & 0.04 & 1.5\% & 0.0\% & 1 & 0.04 & -- & 1.5\% & 0.0\% & 1 & 0.04 & -- \\
a4-32 & 485.5 & 1.9\% & 0.0\% & 1 & 0.15 & 1.6\% & 0.0\% & 1 & 0.08 & -- & 1.6\% & 0.0\% & 1 & 0.08 & -- \\
a4-40 & 557.7 & 1.5\% & 0.0\% & 1 & 0.27 & 1.5\% & 0.0\% & 1 & 0.18 & -- & 1.5\% & 0.0\% & 1 & 0.18 & -- \\
a4-48 & 668.8 & 3.6\% & 0.0\% & 1 & 0.29 & 2.8\% & 0.0\% & 1 & 0.28 & -- & 2.6\% & 0.0\% & 1 & 0.38 & -- \\
a5-40 & 498.4 & 0.6\% & 0.0\% & 1 & 0.14 & 0.5\% & 0.0\% & 1 & 0.15 & -- & 0.5\% & 0.0\% & 1 & 0.15 & -- \\
a5-50 & 686.6 & 4.0\% & 0.3\% & 32 & 1.65 & 3.6\% & 0.0\% & 1 & 0.50 & -70\% & 3.5\% & 0.0\% & 1 & 0.44 & -73\% \\
a5-60 & 808.3 & 2.0\% & 0.0\% & 1 & 0.60 & 1.8\% & 0.0\% & 1 & 0.31 & -48\% & 1.8\% & 0.0\% & 1 & 0.31 & -48\% \\
a6-48 & 604.1 & 1.5\% & 0.0\% & 1 & 0.33 & 1.5\% & 0.0\% & 1 & 0.30 & -- & 1.4\% & 0.0\% & 1 & 0.30 & -- \\
a6-60 & 819.3 & 3.0\% & 0.0\% & 1 & 1.64 & 2.5\% & 0.0\% & 1 & 0.47 & -71\% & 2.5\% & 0.0\% & 1 & 0.63 & -62\% \\
a6-72 & 916.1 & 3.0\% & 0.9\% & 457 & 2.60 & 2.6\% & 0.0\% & 1 & 1.60 & -38\% & 2.5\% & 0.0\% & 1 & 1.30 & -50\% \\
a7-56 & 724.0 & 1.7\% & 0.1\% & 5 & 1.86 & 1.4\% & 0.0\% & 1 & 0.63 & -66\% & 1.4\% & 0.0\% & 1 & 0.63 & -66\% \\
a7-70 & 875.7 & 2.1\% & 0.0\% & 1 & 1.19 & 2.0\% & 0.0\% & 1 & 0.91 & -24\% & 1.9\% & 0.0\% & 1 & 1.07 & -10\% \\
a7-84 & 1033.3 & 2.6\% & 0.1\% & 27 & 2.07 & 2.5\% & 0.0\% & 1 & 1.51 & -27\% & 2.5\% & 0.0\% & 1 & 1.20 & -42\% \\
a8-64 & 747.5 & 2.0\% & 0.0\% & 1 & 3.32 & 1.9\% & 0.0\% & 1 & 0.72 & -78\% & 1.9\% & 0.0\% & 1 & 0.69 & -79\% \\
a8-80 & 945.8 & 3.6\% & 0.6\% & 224 & 3.57 & 3.3\% & 0.0\% & 1 & 1.46 & -59\% & 3.3\% & 0.0\% & 1 & 1.33 & -63\% \\
a8-96 & 1229.7 & 3.7\% & 1.0\% & 1112 & 16.17 & 3.4\% & 0.0\% & 1 & 4.09 & -75\% & 3.4\% & 0.0\% & 24 & 3.97 & -75\% \\
\midrule
Total &  & 2.1\% & 0.1\% & 1875 & 36.5 & 1.9\% & 0.0\% & 24 & 13.7 & -62\% & 1.9\% & 0.0\% & 47 & 13.1 & -64\% \\
\bottomrule
\end{tabular}
\par\vspace{3pt}
\begin{minipage}{\textwidth}\scriptsize
Notes: means over three seeds of LPgap, rGap, Nodes and CPU. Obj.: optimal value, or the best incumbent over all configurations and seeds if no configuration proved optimality. LPgap: gap of the LP bound to the objective value of the run. rGap: gap of the bound when Gurobi leaves the root node to branch, 0 if the run is solved at the root node; --: objective zero or root node not processed. Nodes: branch-and-bound nodes (1: solved at the root node, 0: solved in presolve without a root node). CPU: pure Gurobi time in seconds (graph and model building excluded); TL (x\%) = time limit of 7200~s reached in every seed, largest final gap. $\Delta$CPU: percentage change of CPU against base; --: both below 0.5~s or both at the limit. Total: sums of CPU and Nodes, means of the gaps, percentage change of the summed CPU.
\end{minipage}
\end{table}

\begin{table}[htbp]
\centering\scriptsize\setlength{\tabcolsep}{1.2pt}\renewcommand{\arraystretch}{0.9}
\caption{Base versus base + cuts on the b instances (routing costs~\eqref{eq:laeb:obj}).}
\label{tab:ridecut_LAEB_b_study}
\begin{tabular}{l r rrrr rrrrr rrrrr}
\toprule
 & & \multicolumn{4}{c}{base} & \multicolumn{5}{c}{base + \eqref{eq:ridecut}} & \multicolumn{5}{c}{base + \eqref{eq:lridecut}} \\
\cmidrule(lr){3-6}\cmidrule(lr){7-11}\cmidrule(lr){12-16}
Inst. & Obj. & LPgap & rGap & Nodes & CPU & LPgap & rGap & Nodes & CPU & $\Delta$CPU & LPgap & rGap & Nodes & CPU & $\Delta$CPU \\
\midrule
b2-16 & 309.4 & 0.6\% & 0.0\% & 1 & 0.06 & 0.6\% & 0.0\% & 1 & 0.06 & -- & 0.6\% & 0.0\% & 1 & 0.05 & -- \\
b2-20 & 332.7 & 0.0\% & 0.0\% & 1 & 0.02 & 0.0\% & 0.0\% & 1 & 0.02 & -- & 0.0\% & 0.0\% & 1 & 0.01 & -- \\
b2-24 & 444.7 & 1.3\% & 0.0\% & 1 & 0.07 & 1.3\% & 0.0\% & 1 & 0.12 & -- & 1.3\% & 0.0\% & 1 & 0.06 & -- \\
b3-18 & 301.6 & 0.9\% & 0.0\% & 1 & 0.07 & 0.9\% & 0.0\% & 1 & 0.05 & -- & 0.9\% & 0.0\% & 1 & 0.05 & -- \\
b3-24 & 394.5 & 1.9\% & 0.0\% & 1 & 0.16 & 1.9\% & 0.0\% & 1 & 0.14 & -- & 1.9\% & 0.0\% & 1 & 0.12 & -- \\
b3-30 & 531.4 & 0.1\% & 0.0\% & 1 & 0.04 & 0.1\% & 0.0\% & 1 & 0.05 & -- & 0.1\% & 0.0\% & 1 & 0.05 & -- \\
b3-36 & 603.8 & 0.1\% & 0.0\% & 1 & 0.06 & 0.1\% & 0.0\% & 1 & 0.06 & -- & 0.1\% & 0.0\% & 1 & 0.07 & -- \\
b4-16 & 296.9 & 0.6\% & 0.0\% & 1 & 0.02 & 0.6\% & 0.0\% & 1 & 0.05 & -- & 0.6\% & 0.0\% & 1 & 0.02 & -- \\
b4-24 & 371.4 & 0.1\% & 0.0\% & 1 & 0.08 & 0.1\% & 0.0\% & 1 & 0.07 & -- & 0.1\% & 0.0\% & 1 & 0.05 & -- \\
b4-32 & 494.9 & 0.0\% & 0.0\% & 1 & 0.04 & 0.0\% & 0.0\% & 1 & 0.06 & -- & 0.0\% & 0.0\% & 1 & 0.04 & -- \\
b4-40 & 656.6 & 0.3\% & 0.0\% & 1 & 0.11 & 0.3\% & 0.0\% & 1 & 0.10 & -- & 0.3\% & 0.0\% & 1 & 0.11 & -- \\
b4-48 & 673.8 & 0.7\% & 0.0\% & 1 & 0.93 & 0.7\% & 0.0\% & 1 & 0.70 & -25\% & 0.7\% & 0.0\% & 1 & 0.49 & -47\% \\
b5-40 & 613.7 & 0.4\% & 0.0\% & 1 & 0.14 & 0.4\% & 0.0\% & 1 & 0.17 & -- & 0.4\% & 0.0\% & 1 & 0.17 & -- \\
b5-50 & 761.4 & 1.0\% & 0.0\% & 1 & 0.33 & 1.0\% & 0.0\% & 1 & 0.29 & -- & 1.0\% & 0.0\% & 1 & 0.35 & -- \\
b5-60 & 902.0 & 1.8\% & 0.0\% & 1 & 1.21 & 1.8\% & 0.0\% & 1 & 0.81 & -33\% & 1.8\% & 0.0\% & 1 & 0.98 & -19\% \\
b6-48 & 714.8 & 0.4\% & 0.0\% & 1 & 0.20 & 0.3\% & 0.0\% & 1 & 0.19 & -- & 0.3\% & 0.0\% & 1 & 0.18 & -- \\
b6-60 & 860.0 & 0.1\% & 0.0\% & 1 & 0.26 & 0.1\% & 0.0\% & 1 & 0.26 & -- & 0.1\% & 0.0\% & 1 & 0.23 & -- \\
b6-72 & 978.5 & 0.9\% & 0.1\% & 72 & 3.09 & 0.9\% & 0.0\% & 1 & 2.94 & -5\% & 0.9\% & 0.0\% & 1 & 2.47 & -20\% \\
b7-56 & 824.0 & 0.8\% & 0.0\% & 1 & 1.76 & 0.8\% & 0.0\% & 1 & 1.86 & +6\% & 0.8\% & 0.0\% & 1 & 2.77 & +57\% \\
b7-70 & 912.6 & 1.8\% & 0.0\% & 1 & 0.81 & 1.8\% & 0.0\% & 1 & 0.86 & +6\% & 1.7\% & 0.0\% & 1 & 0.96 & +19\% \\
b7-84 & 1203.3 & 0.8\% & 0.0\% & 1 & 1.19 & 0.8\% & 0.0\% & 1 & 1.24 & +4\% & 0.8\% & 0.0\% & 1 & 1.84 & +55\% \\
b8-64 & 839.9 & 1.6\% & 0.0\% & 1 & 1.00 & 1.6\% & 0.0\% & 1 & 0.69 & -31\% & 1.5\% & 0.0\% & 1 & 0.67 & -33\% \\
b8-80 & 1036.4 & 0.8\% & 0.0\% & 1 & 1.10 & 0.8\% & 0.0\% & 1 & 0.59 & -46\% & 0.8\% & 0.0\% & 1 & 0.69 & -37\% \\
b8-96 & 1185.6 & 1.2\% & 0.0\% & 1 & 2.71 & 1.2\% & 0.0\% & 1 & 2.45 & -10\% & 1.1\% & 0.0\% & 1 & 2.52 & -7\% \\
\midrule
Total &  & 0.8\% & 0.0\% & 95 & 15.5 & 0.7\% & 0.0\% & 24 & 13.8 & -11\% & 0.7\% & 0.0\% & 24 & 14.9 & -4\% \\
\bottomrule
\end{tabular}
\par\vspace{3pt}
\begin{minipage}{\textwidth}\scriptsize
Notes: means over three seeds of LPgap, rGap, Nodes and CPU. Obj.: optimal value, or the best incumbent over all configurations and seeds if no configuration proved optimality. LPgap: gap of the LP bound to the objective value of the run. rGap: gap of the bound when Gurobi leaves the root node to branch, 0 if the run is solved at the root node; --: objective zero or root node not processed. Nodes: branch-and-bound nodes (1: solved at the root node, 0: solved in presolve without a root node). CPU: pure Gurobi time in seconds (graph and model building excluded); TL (x\%) = time limit of 7200~s reached in every seed, largest final gap. $\Delta$CPU: percentage change of CPU against base; --: both below 0.5~s or both at the limit. Total: sums of CPU and Nodes, means of the gaps, percentage change of the summed CPU.
\end{minipage}
\end{table}

\begin{table}[htbp]
\centering\scriptsize\setlength{\tabcolsep}{1.2pt}\renewcommand{\arraystretch}{0.9}
\caption{Base versus base + cuts on the a-X$^\star$ instances (routing costs~\eqref{eq:laeb:obj}).}
\label{tab:ridecut_LAEB_aX_study}
\begin{tabular}{l r rrrr rrrrr rrrrr}
\toprule
 & & \multicolumn{4}{c}{base} & \multicolumn{5}{c}{base + \eqref{eq:ridecut}} & \multicolumn{5}{c}{base + \eqref{eq:lridecut}} \\
\cmidrule(lr){3-6}\cmidrule(lr){7-11}\cmidrule(lr){12-16}
Inst. & Obj. & LPgap & rGap & Nodes & CPU & LPgap & rGap & Nodes & CPU & $\Delta$CPU & LPgap & rGap & Nodes & CPU & $\Delta$CPU \\
\midrule
a2-16 & 278.2 & 0.1\% & 0.0\% & 1 & 0.03 & 0.1\% & 0.0\% & 1 & 0.08 & -- & 0.1\% & 0.0\% & 1 & 0.02 & -- \\
a2-20 & 330.7 & 0.8\% & 0.0\% & 1 & 0.05 & 0.8\% & 0.0\% & 1 & 0.04 & -- & 0.8\% & 0.0\% & 1 & 0.04 & -- \\
a2-24 & 389.1 & 2.2\% & 0.0\% & 1 & 0.18 & 2.0\% & 0.0\% & 1 & 0.15 & -- & 2.0\% & 0.0\% & 1 & 0.14 & -- \\
a3-18 & 272.7 & 6.7\% & 0.0\% & 1 & 0.12 & 5.7\% & 0.0\% & 1 & 0.13 & -- & 5.2\% & 0.0\% & 1 & 0.09 & -- \\
a3-24 & 289.6 & 3.3\% & 1.3\% & 123 & 0.59 & 2.6\% & 0.0\% & 1 & 0.53 & -10\% & 2.6\% & 0.0\% & 1 & 0.46 & -22\% \\
a3-30 & 452.8 & 5.3\% & 0.0\% & 1 & 0.38 & 4.3\% & 0.0\% & 1 & 0.11 & -- & 4.3\% & 0.0\% & 1 & 0.12 & -- \\
a3-36 & 501.0 & 1.5\% & 0.2\% & 50 & 0.44 & 1.3\% & 0.0\% & 1 & 0.21 & -- & 1.3\% & 0.0\% & 1 & 0.20 & -- \\
a4-16 & 235.2 & 5.8\% & 0.0\% & 1 & 0.39 & 5.6\% & 0.0\% & 1 & 0.21 & -- & 5.6\% & 0.0\% & 1 & 0.20 & -- \\
a4-24 & 359.4 & 3.6\% & 0.0\% & 1 & 0.17 & 2.5\% & 0.0\% & 1 & 0.12 & -- & 2.5\% & 0.0\% & 1 & 0.12 & -- \\
a4-32 & 447.3 & 3.7\% & 1.6\% & 208 & 1.16 & 3.2\% & 0.0\% & 1 & 0.48 & -59\% & 3.2\% & 0.0\% & 1 & 0.33 & -72\% \\
a4-40 & 509.0 & 2.4\% & 0.4\% & 119 & 2.09 & 2.0\% & 0.0\% & 1 & 1.24 & -41\% & 2.0\% & 0.0\% & 1 & 0.91 & -56\% \\
a4-48 & 620.3 & 7.8\% & 3.8\% & 2052 & 21.39 & 6.2\% & 0.2\% & 24 & 5.09 & -76\% & 6.2\% & 0.5\% & 48 & 7.07 & -67\% \\
a5-40 & 464.0 & 3.0\% & 0.7\% & 174 & 1.65 & 2.7\% & 0.0\% & 1 & 0.75 & -55\% & 2.4\% & 0.0\% & 1 & 0.81 & -51\% \\
a5-50 & 621.9 & 5.2\% & 2.6\% & 1711 & 16.45 & 4.4\% & 0.7\% & 109 & 2.63 & -84\% & 4.4\% & 0.0\% & 1 & 3.00 & -82\% \\
a5-60 & 745.4 & 4.7\% & 1.3\% & 431 & 4.06 & 4.2\% & 0.0\% & 1 & 2.48 & -39\% & 4.0\% & 0.0\% & 1 & 2.70 & -33\% \\
a6-48 & 572.5 & 5.4\% & 2.4\% & 2034 & 73.69 & 4.8\% & 0.3\% & 48 & 7.20 & -90\% & 4.6\% & 0.5\% & 147 & 7.06 & -90\% \\
a6-60 & 757.9 & 3.9\% & 0.9\% & 1892 & 24.92 & 3.5\% & 0.0\% & 1 & 8.74 & -65\% & 3.1\% & 0.1\% & 30 & 5.73 & -77\% \\
a6-72 & 868.4 & 6.2\% & 3.5\% & 33716 & 1756.70 & 5.1\% & 1.6\% & 731 & 40.98 & -98\% & 5.1\% & 1.6\% & 820 & 45.02 & -97\% \\
a7-56 & 663.5 & 4.7\% & 2.4\% & 2574 & 90.39 & 4.1\% & 1.1\% & 385 & 10.51 & -88\% & 4.0\% & 1.2\% & 316 & 9.11 & -90\% \\
a7-70 & 815.3 & 6.7\% & 1.5\% & 1075 & 123.83 & 5.8\% & 0.0\% & 1 & 27.48 & -78\% & 5.6\% & 0.0\% & 1 & 24.56 & -80\% \\
a7-84 & 950.6 & 7.6\% & 3.3\% & 26914 & 3844.55 & 6.1\% & 0.7\% & 961 & 81.64 & -98\% & 5.9\% & 0.5\% & 1178 & 84.05 & -98\% \\
a8-64 & 701.2 & 4.8\% & 1.4\% & 2507 & 261.45 & 4.2\% & 0.7\% & 572 & 48.51 & -81\% & 4.1\% & 0.5\% & 560 & 48.22 & -82\% \\
a8-80 & 880.3 & 7.3\% & 3.6\% & 50369 & 5159.08 & 5.7\% & 0.6\% & 596 & 73.50 & -99\% & 5.5\% & 0.7\% & 516 & 73.38 & -99\% \\
a8-96 & 1115.1 & 5.8\% & 3.8\% & 22327 & TL (1.7\%) & 5.4\% & 1.5\% & 2329 & 490.44 & -93\% & 5.4\% & 1.4\% & 1765 & 442.12 & -94\% \\
\midrule
Total &  & 4.5\% & 1.4\% & 148283 & 18583.8 & 3.8\% & 0.3\% & 5770 & 803.3 & -96\% & 3.8\% & 0.3\% & 5395 & 755.4 & -96\% \\
\bottomrule
\end{tabular}
\par\vspace{3pt}
\begin{minipage}{\textwidth}\scriptsize
Notes: means over three seeds of LPgap, rGap, Nodes and CPU. Obj.: optimal value, or the best incumbent over all configurations and seeds if no configuration proved optimality. LPgap: gap of the LP bound to the objective value of the run. rGap: gap of the bound when Gurobi leaves the root node to branch, 0 if the run is solved at the root node; --: objective zero or root node not processed. Nodes: branch-and-bound nodes (1: solved at the root node, 0: solved in presolve without a root node). CPU: pure Gurobi time in seconds (graph and model building excluded); TL (x\%) = time limit of 7200~s reached in every seed, largest final gap. $\Delta$CPU: percentage change of CPU against base; --: both below 0.5~s or both at the limit. Total: sums of CPU and Nodes, means of the gaps, percentage change of the summed CPU.
\end{minipage}
\end{table}

\begin{table}[htbp]
\centering\scriptsize\setlength{\tabcolsep}{1.2pt}\renewcommand{\arraystretch}{0.9}
\caption{Base versus base + cuts on the b-X$^\star$ instances (routing costs~\eqref{eq:laeb:obj}).}
\label{tab:ridecut_LAEB_bX_study}
\begin{tabular}{l r rrrr rrrrr rrrrr}
\toprule
 & & \multicolumn{4}{c}{base} & \multicolumn{5}{c}{base + \eqref{eq:ridecut}} & \multicolumn{5}{c}{base + \eqref{eq:lridecut}} \\
\cmidrule(lr){3-6}\cmidrule(lr){7-11}\cmidrule(lr){12-16}
Inst. & Obj. & LPgap & rGap & Nodes & CPU & LPgap & rGap & Nodes & CPU & $\Delta$CPU & LPgap & rGap & Nodes & CPU & $\Delta$CPU \\
\midrule
b2-16 & 282.5 & 1.5\% & 0.0\% & 1 & 0.26 & 1.5\% & 0.0\% & 1 & 0.22 & -- & 1.5\% & 0.0\% & 1 & 0.17 & -- \\
b2-20 & 323.6 & 0.8\% & 0.0\% & 1 & 0.07 & 0.8\% & 0.0\% & 1 & 0.04 & -- & 0.8\% & 0.0\% & 1 & 0.04 & -- \\
b2-24 & 412.3 & 0.0\% & 0.0\% & 1 & 0.06 & 0.0\% & 0.0\% & 1 & 0.06 & -- & 0.0\% & 0.0\% & 1 & 0.05 & -- \\
b3-18 & 290.4 & 1.5\% & 0.0\% & 1 & 0.07 & 1.5\% & 0.0\% & 1 & 0.09 & -- & 1.5\% & 0.0\% & 1 & 0.07 & -- \\
b3-24 & 363.7 & 0.5\% & 0.0\% & 1 & 0.11 & 0.5\% & 0.0\% & 1 & 0.19 & -- & 0.5\% & 0.0\% & 1 & 0.17 & -- \\
b3-30 & 504.3 & 1.8\% & 0.0\% & 1 & 0.38 & 1.8\% & 0.0\% & 1 & 0.28 & -- & 1.8\% & 0.0\% & 1 & 0.22 & -- \\
b3-36 & 565.9 & 0.0\% & 0.0\% & 1 & 0.07 & 0.0\% & 0.0\% & 1 & 0.08 & -- & 0.0\% & 0.0\% & 1 & 0.08 & -- \\
b4-16 & 289.9 & 1.0\% & 0.0\% & 1 & 0.03 & 1.0\% & 0.0\% & 1 & 0.04 & -- & 1.0\% & 0.0\% & 1 & 0.02 & -- \\
b4-24 & 347.0 & 2.0\% & 0.0\% & 1 & 0.42 & 2.0\% & 0.0\% & 1 & 0.37 & -- & 2.0\% & 0.0\% & 1 & 0.46 & -- \\
b4-32 & 491.0 & 0.7\% & 0.0\% & 1 & 0.17 & 0.7\% & 0.0\% & 1 & 0.13 & -- & 0.7\% & 0.0\% & 1 & 0.10 & -- \\
b4-40 & 628.3 & 0.5\% & 0.0\% & 1 & 0.33 & 0.5\% & 0.0\% & 1 & 0.23 & -- & 0.5\% & 0.0\% & 1 & 0.21 & -- \\
b4-48 & 627.4 & 0.9\% & 0.0\% & 1 & 1.23 & 0.9\% & 0.0\% & 1 & 0.98 & -20\% & 0.9\% & 0.0\% & 1 & 1.25 & +2\% \\
b5-40 & 585.1 & 1.8\% & 0.0\% & 1 & 2.22 & 1.7\% & 0.0\% & 1 & 1.58 & -29\% & 1.6\% & 0.0\% & 1 & 2.13 & -4\% \\
b5-50 & 708.8 & 0.6\% & 0.0\% & 1 & 0.49 & 0.6\% & 0.0\% & 1 & 0.92 & +88\% & 0.6\% & 0.0\% & 1 & 0.77 & +57\% \\
b5-60 & 851.9 & 1.4\% & 0.0\% & 1 & 1.70 & 1.4\% & 0.0\% & 1 & 1.59 & -6\% & 1.4\% & 0.0\% & 1 & 1.52 & -11\% \\
b6-48 & 691.6 & 0.9\% & 0.0\% & 1 & 0.38 & 0.9\% & 0.0\% & 1 & 0.41 & -- & 0.9\% & 0.0\% & 1 & 0.44 & -- \\
b6-60 & 841.6 & 1.3\% & 0.0\% & 1 & 0.89 & 1.3\% & 0.0\% & 1 & 0.87 & -2\% & 1.3\% & 0.0\% & 1 & 0.98 & +10\% \\
b6-72 & 930.3 & 1.6\% & 0.0\% & 1 & 4.86 & 1.6\% & 0.0\% & 1 & 6.30 & +30\% & 1.5\% & 0.0\% & 1 & 3.62 & -26\% \\
b7-56 & 787.9 & 0.8\% & 0.0\% & 1 & 4.53 & 0.8\% & 0.0\% & 1 & 5.94 & +31\% & 0.7\% & 0.0\% & 1 & 6.06 & +34\% \\
b7-70 & 865.3 & 1.8\% & 0.0\% & 1 & 1.66 & 1.8\% & 0.1\% & 34 & 2.93 & +77\% & 1.7\% & 0.1\% & 19 & 4.06 & +145\% \\
b7-84 & 1141.2 & 2.4\% & 0.7\% & 1315 & 19.70 & 2.3\% & 0.3\% & 185 & 13.82 & -30\% & 2.2\% & 0.4\% & 320 & 14.19 & -28\% \\
b8-64 & 818.3 & 1.4\% & 0.0\% & 1 & 1.65 & 1.4\% & 0.0\% & 1 & 2.03 & +23\% & 1.4\% & 0.0\% & 1 & 1.99 & +21\% \\
b8-80 & 998.3 & 1.7\% & 0.0\% & 1 & 2.74 & 1.7\% & 0.0\% & 1 & 2.06 & -25\% & 1.7\% & 0.0\% & 1 & 2.11 & -23\% \\
b8-96 & 1137.7 & 2.3\% & 0.1\% & 69 & 42.94 & 2.2\% & 0.0\% & 1 & 34.44 & -20\% & 2.2\% & 0.0\% & 1 & 35.36 & -18\% \\
\midrule
Total &  & 1.2\% & 0.0\% & 1406 & 87.0 & 1.2\% & 0.0\% & 241 & 75.6 & -13\% & 1.2\% & 0.0\% & 361 & 76.1 & -13\% \\
\bottomrule
\end{tabular}
\par\vspace{3pt}
\begin{minipage}{\textwidth}\scriptsize
Notes: means over three seeds of LPgap, rGap, Nodes and CPU. Obj.: optimal value, or the best incumbent over all configurations and seeds if no configuration proved optimality. LPgap: gap of the LP bound to the objective value of the run. rGap: gap of the bound when Gurobi leaves the root node to branch, 0 if the run is solved at the root node; --: objective zero or root node not processed. Nodes: branch-and-bound nodes (1: solved at the root node, 0: solved in presolve without a root node). CPU: pure Gurobi time in seconds (graph and model building excluded); TL (x\%) = time limit of 7200~s reached in every seed, largest final gap. $\Delta$CPU: percentage change of CPU against base; --: both below 0.5~s or both at the limit. Total: sums of CPU and Nodes, means of the gaps, percentage change of the summed CPU.
\end{minipage}
\end{table}

The inequalities act where branching is needed.
On the a instances the base configuration needs 36.5 seconds and 1875 nodes in total, the configuration with the inequalities \eqref{eq:ridecut} 13.7 seconds and 24 nodes, that is, every instance is solved at the root node; the CPU time falls on all ten instances that take more than half a second, by 56~percent on average.
On the a-X$^\star$ instances the base configuration needs 18\,584 seconds and 148\,283 nodes and does not solve a8-96-X$^\star$ within the time limit in any of the three seeds (final gaps 1.5 to 1.7~percent), whereas the configuration with the inequalities \eqref{eq:ridecut} solves every instance, in 803 seconds and 5770 nodes in total; it is faster on all sixteen instances above half a second, by 72~percent on average and by 96~percent in total.
The LP relaxation gap improves only a little, from 4.5 to 3.8~percent on a-X$^\star$, because the time variables do not enter the objective; the effect shows in the root gap after Gurobi's own cutting planes, which falls from 1.4 to 0.3~percent, and in the size of the search trees.
The optimum of a8-96-X$^\star$, 1115.13, is proven in 391 to 564 seconds by the three seeds with the inequalities; it is also proven by GP with the inequalities in 242 seconds and by the plain model with the inequalities in 618 seconds.
The b and b-X$^\star$ instances are solved at or near the root node by the base configuration already (95 and 1406 nodes in total), and the inequalities \eqref{eq:ridecut} change the totals by $-11$ and $-13$~percent with differences of a few tenths of a second per instance in both directions.
The weaker configurations, run once each, show the same picture: on a-X$^\star$ the plain model needs 24\,204 seconds with two unsolved instances and the plain model with the inequalities 902 seconds with none, GP needs 12\,866 seconds with one unsolved instance and GP with the inequalities 464 seconds with none.
Every run that reached the time limit belongs to a configuration without the inequalities.

The strengthened inequalities~\eqref{eq:lridecut} behave like~\eqref{eq:ridecut}, as the small share of arcs with forced waiting suggests.
Over the 144 instances with the routing costs, the LP gaps of the two configurations coincide to one decimal on 110 and are lower with~\eqref{eq:lridecut} on 34, by at most 0.8 percentage points (a4-16-X$^{\star\star}$, 5.4 against 4.6~percent); the root gaps differ by a few tenths in both directions.
The totals move within the variation between the seeds: on the a instances from 13.7 to 13.1 seconds, on a-X$^\star$ from 803 to 755 seconds and from 5770 to 5395 nodes, where~\eqref{eq:lridecut} is faster on nine of the fifteen instances above half a second and proves the optimum of a8-96-X$^\star$ in 416 to 457 seconds, and on b and b-X$^\star$ from 13.8 to 14.9 and from 75.6 to 76.1 seconds.
The weaker configurations with~\eqref{eq:lridecut} solve every a-X$^\star$ instance as well, in 1018 (plain) and 497 seconds (GP) against 902 and 464 seconds with~\eqref{eq:ridecut}.

\paragraph{Wider time windows.}
Tables~\ref{tab:ridecut_LAEB_aX2_study_x2} and~\ref{tab:ridecut_LAEB_bX2_study_x2} repeat the comparison of the base configuration with and without the inequalities on the sets a-X$^{\star\star}$ and b-X$^{\star\star}$.

\begin{table}[htbp]
\centering\scriptsize\setlength{\tabcolsep}{1.2pt}\renewcommand{\arraystretch}{0.9}
\caption{Base versus base + cuts on the a-X$^{\star\star}$ instances (routing costs~\eqref{eq:laeb:obj}).}
\label{tab:ridecut_LAEB_aX2_study_x2}
\begin{tabular}{l r rrrr rrrrr rrrrr}
\toprule
 & & \multicolumn{4}{c}{base} & \multicolumn{5}{c}{base + \eqref{eq:ridecut}} & \multicolumn{5}{c}{base + \eqref{eq:lridecut}} \\
\cmidrule(lr){3-6}\cmidrule(lr){7-11}\cmidrule(lr){12-16}
Inst. & Obj. & LPgap & rGap & Nodes & CPU & LPgap & rGap & Nodes & CPU & $\Delta$CPU & LPgap & rGap & Nodes & CPU & $\Delta$CPU \\
\midrule
a2-16 & 260.6 & 3.1\% & 0.0\% & 1 & 0.04 & 3.1\% & 0.0\% & 1 & 0.03 & -- & 3.1\% & 0.0\% & 1 & 0.03 & -- \\
a2-20 & 317.7 & 5.8\% & 0.0\% & 1 & 0.42 & 5.4\% & 0.0\% & 1 & 0.12 & -- & 5.4\% & 0.0\% & 1 & 0.10 & -- \\
a2-24 & 378.2 & 3.9\% & 0.9\% & 64 & 0.97 & 3.9\% & 0.0\% & 1 & 0.31 & -68\% & 3.9\% & 0.0\% & 1 & 0.38 & -61\% \\
a3-18 & 258.0 & 6.9\% & 3.2\% & 423 & 0.70 & 6.2\% & 0.0\% & 1 & 0.51 & -27\% & 6.2\% & 0.0\% & 1 & 0.50 & -29\% \\
a3-24 & 282.2 & 3.0\% & 1.8\% & 161 & 0.94 & 2.5\% & 0.0\% & 1 & 0.45 & -52\% & 2.5\% & 0.0\% & 1 & 0.44 & -53\% \\
a3-30 & 437.2 & 5.4\% & 3.2\% & 940 & 2.10 & 4.7\% & 0.0\% & 1 & 0.58 & -72\% & 4.7\% & 0.0\% & 1 & 0.71 & -66\% \\
a3-36 & 473.0 & 1.6\% & 0.0\% & 1 & 0.53 & 1.3\% & 0.0\% & 1 & 0.43 & -19\% & 1.3\% & 0.0\% & 1 & 0.39 & -26\% \\
a4-16 & 218.5 & 6.6\% & 4.3\% & 306 & 0.96 & 5.4\% & 0.0\% & 1 & 0.73 & -24\% & 4.6\% & 0.0\% & 1 & 0.70 & -27\% \\
a4-24 & 352.0 & 6.3\% & 1.1\% & 105 & 0.81 & 6.0\% & 0.0\% & 1 & 0.50 & -38\% & 6.0\% & 0.0\% & 1 & 0.54 & -33\% \\
a4-32 & 420.5 & 6.3\% & 4.0\% & 1584 & 13.98 & 5.6\% & 0.3\% & 27 & 4.13 & -70\% & 5.6\% & 0.6\% & 146 & 3.32 & -76\% \\
a4-40 & 492.1 & 5.4\% & 2.0\% & 4908 & 189.85 & 4.6\% & 1.5\% & 971 & 46.29 & -76\% & 4.6\% & 1.5\% & 967 & 50.08 & -74\% \\
a4-48 & 586.9 & 10.4\% & 6.4\% & 66229 & 2230.51 & 8.3\% & 2.1\% & 1105 & 204.39 & -91\% & 8.2\% & 1.7\% & 1213 & 97.30 & -96\% \\
a5-40 & 433.4 & 7.2\% & 3.8\% & 1987 & 44.16 & 6.0\% & 0.3\% & 36 & 3.70 & -92\% & 5.9\% & 0.4\% & 38 & 4.41 & -90\% \\
a5-50 & 605.7 & 8.1\% & 5.6\% & 55922 & 2047.86 & 6.6\% & 2.2\% & 924 & 75.77 & -96\% & 6.6\% & 2.3\% & 1189 & 53.47 & -97\% \\
a5-60 & 696.8 & 9.4\% & 5.2\% & 23817 & 1205.64 & 7.3\% & 1.6\% & 1369 & 91.76 & -92\% & 7.3\% & 1.1\% & 991 & 56.25 & -95\% \\
a6-48 & 535.2 & 6.0\% & 2.0\% & 3036 & 763.01 & 5.1\% & 0.2\% & 142 & 59.44 & -92\% & 4.9\% & 0.0\% & 1 & 24.98 & -97\% \\
a6-60 & 701.6 & 6.9\% & 4.0\% & 27644 & TL (1.3\%) & 5.4\% & 2.2\% & 2959 & 737.69 & -90\% & 5.2\% & 2.1\% & 1526 & 760.16 & -89\% \\
a6-72 & 797.8 & 7.6\% & 5.3\% & 24262 & TL (2.3\%) & 5.7\% & 1.9\% & 815 & 309.15 & -96\% & 5.6\% & 1.8\% & 1566 & 387.09 & -95\% \\
a7-56 & 626.5 & 4.8\% & 1.8\% & 2198 & 459.97 & 4.3\% & 0.0\% & 1 & 33.48 & -93\% & 4.2\% & 0.0\% & 1 & 22.07 & -95\% \\
a7-70 & 777.9 & 9.5\% & 6.3\% & 12344 & TL (3.4\%) & 7.9\% & 2.9\% & 7054 & 2581.81 & -64\% & 7.8\% & 3.0\% & 5149 & 1853.28 & -74\% \\
a7-84 & 899.3 & 12.0\% & 8.7\% & 9230 & TL (7.1\%) & 8.4\% & 3.5\% & 12949 & 5382.39 & -25\% & 8.3\% & 3.6\% & 15982 & 5544.69 & -23\% \\
a8-64 & 670.0 & 6.9\% & 4.8\% & 14639 & TL (1.5\%) & 5.5\% & 2.0\% & 990 & 1122.42 & -84\% & 5.5\% & 1.5\% & 1213 & 569.91 & -92\% \\
a8-80 & 834.3 & 11.2\% & 7.3\% & 7982 & TL (5.8\%) & 9.1\% & 3.9\% & 10202 & TL (0.8\%) & -- & 9.1\% & 4.0\% & 13643 & TL (0.8\%) & -- \\
a8-96 & 1046.0 & 12.5\% & 9.0\% & 2372 & TL (8.2\%) & 8.3\% & 3.3\% & 5489 & TL (0.8\%) & -- & 8.2\% & 3.3\% & 7628 & TL (0.7\%) & -- \\
\midrule
Total &  & 7.0\% & 3.8\% & 260156 & 57362.5 & 5.7\% & 1.2\% & 45042 & 25056.1 & -56\% & 5.6\% & 1.1\% & 51262 & 23830.8 & -58\% \\
\bottomrule
\end{tabular}
\par\vspace{3pt}
\begin{minipage}{\textwidth}\scriptsize
Notes: means over three seeds of LPgap, rGap, Nodes and CPU. Obj.: optimal value, or the best incumbent over all configurations and seeds if no configuration proved optimality. LPgap: gap of the LP bound to the objective value of the run. rGap: gap of the bound when Gurobi leaves the root node to branch, 0 if the run is solved at the root node; --: objective zero or root node not processed. Nodes: branch-and-bound nodes (1: solved at the root node, 0: solved in presolve without a root node). CPU: pure Gurobi time in seconds (graph and model building excluded); TL (x\%) = time limit of 7200~s reached in every seed, largest final gap. $\Delta$CPU: percentage change of CPU against base; --: both below 0.5~s or both at the limit. Total: sums of CPU and Nodes, means of the gaps, percentage change of the summed CPU.
\end{minipage}
\end{table}

\begin{table}[htbp]
\centering\scriptsize\setlength{\tabcolsep}{1.2pt}\renewcommand{\arraystretch}{0.9}
\caption{Base versus base + cuts on the b-X$^{\star\star}$ instances (routing costs~\eqref{eq:laeb:obj}).}
\label{tab:ridecut_LAEB_bX2_study_x2}
\begin{tabular}{l r rrrr rrrrr rrrrr}
\toprule
 & & \multicolumn{4}{c}{base} & \multicolumn{5}{c}{base + \eqref{eq:ridecut}} & \multicolumn{5}{c}{base + \eqref{eq:lridecut}} \\
\cmidrule(lr){3-6}\cmidrule(lr){7-11}\cmidrule(lr){12-16}
Inst. & Obj. & LPgap & rGap & Nodes & CPU & LPgap & rGap & Nodes & CPU & $\Delta$CPU & LPgap & rGap & Nodes & CPU & $\Delta$CPU \\
\midrule
b2-16 & 262.3 & 1.8\% & 0.0\% & 1 & 0.22 & 1.8\% & 0.0\% & 1 & 0.17 & -- & 1.8\% & 0.0\% & 1 & 0.24 & -- \\
b2-20 & 301.4 & 0.0\% & 0.0\% & 1 & 0.03 & 0.0\% & 0.0\% & 1 & 0.03 & -- & 0.0\% & 0.0\% & 1 & 0.03 & -- \\
b2-24 & 404.5 & 2.5\% & 0.0\% & 1 & 0.11 & 2.5\% & 0.0\% & 1 & 0.11 & -- & 2.3\% & 0.0\% & 1 & 0.10 & -- \\
b3-18 & 281.7 & 0.9\% & 0.0\% & 1 & 0.54 & 0.9\% & 0.0\% & 1 & 0.38 & -30\% & 0.9\% & 0.0\% & 1 & 0.22 & -59\% \\
b3-24 & 356.0 & 2.0\% & 0.0\% & 1 & 0.43 & 2.0\% & 0.0\% & 1 & 0.50 & +16\% & 2.0\% & 0.0\% & 1 & 0.59 & +37\% \\
b3-30 & 494.1 & 2.7\% & 0.0\% & 1 & 0.33 & 2.7\% & 0.0\% & 1 & 0.59 & +79\% & 2.7\% & 0.0\% & 1 & 0.44 & -- \\
b3-36 & 537.5 & 0.0\% & 0.0\% & 1 & 0.18 & 0.0\% & 0.0\% & 1 & 0.17 & -- & 0.0\% & 0.0\% & 1 & 0.17 & -- \\
b4-16 & 275.4 & 0.3\% & 0.0\% & 1 & 0.04 & 0.3\% & 0.0\% & 1 & 0.04 & -- & 0.3\% & 0.0\% & 1 & 0.04 & -- \\
b4-24 & 332.2 & 0.2\% & 0.0\% & 1 & 0.29 & 0.2\% & 0.0\% & 1 & 0.35 & -- & 0.2\% & 0.0\% & 1 & 0.34 & -- \\
b4-32 & 475.8 & 0.5\% & 0.0\% & 1 & 0.11 & 0.5\% & 0.0\% & 1 & 0.11 & -- & 0.5\% & 0.0\% & 1 & 0.11 & -- \\
b4-40 & 611.1 & 1.1\% & 0.0\% & 1 & 0.87 & 1.1\% & 0.0\% & 1 & 1.01 & +16\% & 1.1\% & 0.0\% & 1 & 0.83 & -5\% \\
b4-48 & 611.8 & 2.6\% & 0.1\% & 18 & 11.96 & 2.3\% & 0.0\% & 1 & 9.74 & -19\% & 2.3\% & 0.0\% & 1 & 13.12 & +10\% \\
b5-40 & 553.8 & 0.1\% & 0.0\% & 1 & 0.85 & 0.1\% & 0.0\% & 1 & 1.27 & +49\% & 0.1\% & 0.0\% & 1 & 1.06 & +25\% \\
b5-50 & 700.2 & 3.2\% & 1.3\% & 1574 & 11.79 & 3.1\% & 0.6\% & 224 & 4.86 & -59\% & 3.1\% & 1.0\% & 366 & 4.58 & -61\% \\
b5-60 & 816.9 & 2.3\% & 0.4\% & 125 & 8.74 & 2.3\% & 0.2\% & 81 & 11.40 & +30\% & 2.3\% & 0.2\% & 49 & 9.88 & +13\% \\
b6-48 & 666.6 & 0.9\% & 0.0\% & 1 & 0.76 & 0.9\% & 0.0\% & 1 & 0.84 & +11\% & 0.9\% & 0.0\% & 1 & 0.77 & +1\% \\
b6-60 & 819.1 & 2.6\% & 0.1\% & 31 & 5.33 & 2.6\% & 0.0\% & 1 & 5.28 & -1\% & 2.6\% & 0.0\% & 1 & 4.91 & -8\% \\
b6-72 & 892.3 & 2.6\% & 0.6\% & 1151 & 36.54 & 2.5\% & 0.4\% & 343 & 37.88 & +4\% & 2.5\% & 0.7\% & 514 & 19.54 & -47\% \\
b7-56 & 766.9 & 1.2\% & 0.0\% & 1 & 19.62 & 1.2\% & 0.1\% & 34 & 32.42 & +65\% & 1.2\% & 0.2\% & 56 & 42.86 & +118\% \\
b7-70 & 833.3 & 2.1\% & 0.2\% & 25 & 11.92 & 2.1\% & 0.1\% & 5 & 8.29 & -30\% & 2.1\% & 0.2\% & 44 & 8.40 & -30\% \\
b7-84 & 1091.8 & 3.6\% & 2.0\% & 28732 & 2491.71 & 3.3\% & 1.6\% & 6247 & 814.18 & -67\% & 3.3\% & 1.6\% & 6190 & 759.97 & -70\% \\
b8-64 & 794.7 & 2.1\% & 0.1\% & 23 & 17.54 & 2.0\% & 0.1\% & 15 & 22.91 & +31\% & 2.0\% & 0.1\% & 20 & 20.38 & +16\% \\
b8-80 & 975.9 & 2.6\% & 0.3\% & 238 & 22.68 & 2.6\% & 0.3\% & 258 & 13.44 & -41\% & 2.6\% & 0.3\% & 231 & 11.98 & -47\% \\
b8-96 & 1091.2 & 2.1\% & 0.8\% & 1516 & 457.20 & 2.0\% & 0.0\% & 1 & 123.69 & -73\% & 2.0\% & 0.0\% & 1 & 96.77 & -79\% \\
\midrule
Total &  & 1.7\% & 0.2\% & 33447 & 3099.8 & 1.6\% & 0.1\% & 7223 & 1089.7 & -65\% & 1.6\% & 0.2\% & 7486 & 997.3 & -68\% \\
\bottomrule
\end{tabular}
\par\vspace{3pt}
\begin{minipage}{\textwidth}\scriptsize
Notes: means over three seeds of LPgap, rGap, Nodes and CPU. Obj.: optimal value, or the best incumbent over all configurations and seeds if no configuration proved optimality. LPgap: gap of the LP bound to the objective value of the run. rGap: gap of the bound when Gurobi leaves the root node to branch, 0 if the run is solved at the root node; --: objective zero or root node not processed. Nodes: branch-and-bound nodes (1: solved at the root node, 0: solved in presolve without a root node). CPU: pure Gurobi time in seconds (graph and model building excluded); TL (x\%) = time limit of 7200~s reached in every seed, largest final gap. $\Delta$CPU: percentage change of CPU against base; --: both below 0.5~s or both at the limit. Total: sums of CPU and Nodes, means of the gaps, percentage change of the summed CPU.
\end{minipage}
\end{table}

With windows of 45 minutes the a instances become considerably harder, and the gap between the two configurations widens.
The base configuration leaves seven of the 24 a-X$^{\star\star}$ instances unsolved in every seed (a6-60, a6-72, a7-70, a7-84, a8-64, a8-80 and a8-96, final gaps 0.9 to 8.2~percent) and finds no feasible solution at all for a8-96-X$^{\star\star}$ in two of the three seeds; the objective value reported for it is that of the remaining seed.
With the inequalities \eqref{eq:ridecut}, five of these seven instances are solved in every seed, in 309 to 5382 seconds on average, and only a8-80-X$^{\star\star}$ and a8-96-X$^{\star\star}$ remain open, with final gaps of 0.7 to 0.8 and 0.6 to 0.8~percent; for a8-96-X$^{\star\star}$ all three seeds end with the incumbent 1046.04, against 1083.30 for the only base run that found a solution.
The configuration with the inequalities is faster on all twenty instances on which at least one configuration solves the instance and takes more than half a second, by 68~percent on average, and needs 25\,056 instead of 57\,363 seconds and 45\,042 instead of 260\,156 nodes in total, where the totals count the time-limit runs with 7200 seconds and thus understate the difference.
As on a-X$^\star$, the LP relaxation gap changes little (7.0 to 5.7~percent on average), whereas the root gap falls from 3.8 to 1.2~percent.
On b-X$^{\star\star}$ both configurations solve every instance; the totals fall from 3100 to 1090 seconds and from 33\,447 to 7223 nodes, mainly through the instances with large search trees (b7-84-X$^{\star\star}$ from 2492 to 814 seconds, b8-96-X$^{\star\star}$ from 457 to 124 seconds and from 1516 nodes to the root node, b5-50-X$^{\star\star}$ from 12 to 5 seconds), while the instances solved at or near the root node show the same small differences in both directions as on b-X$^\star$ (eight faster, nine slower, $-1$~percent on average).
With the strengthened inequalities~\eqref{eq:lridecut} the picture is the same: the same five instances are solved, in 387 to 5545 seconds on average, a8-80-X$^{\star\star}$ and a8-96-X$^{\star\star}$ remain open with the same incumbents and final gaps of 0.3 to 0.8 and 0.6 to 0.7~percent, and the totals are 23\,831 seconds and 51\,262 nodes on a-X$^{\star\star}$ (5~percent less time than with~\eqref{eq:ridecut}, faster on ten of seventeen instances) and 997 seconds and 7486 nodes on b-X$^{\star\star}$ (8~percent less, faster on twelve of sixteen).

\section{Conclusion} \label{sec:conclusion}
The paper presents new valid inequalities for the event-based mixed-integer programming formulation for the Dial-a-Ride Problem. A substantial weakness of this formulation is the implementation of time progression via Big $M$ constraints which need to hold in the extreme case and lead, therefore, often to a weak LP relaxation.

Our introduced valid inequalities \eqref{eq:ridecut} and the improved version \eqref{eq:lridecut}, which additionally include enforced waiting times, address exactly this issue for the relation between a request's pickup and delivery location. For these locations we know that the delivery location is visited after the pickup location within the same tour. Moreover, the event structure allows us to identify all visited events in between. These relations are used for the valid inequalities.

The valid inequalities are evaluated in a comprehensive computational study on the well-known DARP benchmark instances where computation times could be reduced by up to 99\%. Since all $X^\star$ instances could be solved with the inequalities, we extend their time windows by another 15 minutes. Finally, we also applied the valid inequalities to a time-related objective function. Indeed, the LP relaxation bound is zero for many instances and this objective. Nevertheless, the new valid inequalities are able to increase the root node bound and decrease the required search tree size.
The strengthened inequalities \eqref{eq:lridecut} dominate \eqref{eq:ridecut} in theory, but forced waiting is rare after the graph preprocessing, so that the LP bounds of the two versions coincide on most instances and the computation times vary, slightly in favour of \eqref{eq:lridecut} for the routing costs and of \eqref{eq:ridecut} for the relative detour.

However, we only tested the LAEB formulation and instances with up to 96 customers.

Our valid inequalities work in general for the event-based model formulation and can, thereby, also be applied to other problem variants like those named in the introduction. Future research, thus, should use the valid inequalities to improve on model formulations for these problems.

\appendix
\section{Relative detour}
\label{sec:results:detour}

Tables~\ref{tab:ridecut_LAEB_a_study_detour}--\ref{tab:ridecut_LAEB_bX_study_detour} report the same comparison as in Section \ref{sec:results:cost} for the relative-detour objective~\eqref{eq:detour}, again with the same configurations; the objective values are given with two decimals.

\begin{table}[htbp]
\centering\scriptsize\setlength{\tabcolsep}{1.2pt}\renewcommand{\arraystretch}{0.9}
\caption{Base versus base + cuts on the a instances (relative detour~\eqref{eq:detour}).}
\label{tab:ridecut_LAEB_a_study_detour}
\begin{tabular}{l r rrrr rrrrr rrrrr}
\toprule
 & & \multicolumn{4}{c}{base} & \multicolumn{5}{c}{base + \eqref{eq:ridecut}} & \multicolumn{5}{c}{base + \eqref{eq:lridecut}} \\
\cmidrule(lr){3-6}\cmidrule(lr){7-11}\cmidrule(lr){12-16}
Inst. & Obj. & LPgap & rGap & Nodes & CPU & LPgap & rGap & Nodes & CPU & $\Delta$CPU & LPgap & rGap & Nodes & CPU & $\Delta$CPU \\
\midrule
a2-16 & 1.27 & 79.4\% & 0.0\% & 1 & 0.04 & 75.8\% & 0.0\% & 1 & 0.02 & -- & 75.8\% & 0.0\% & 1 & 0.02 & -- \\
a2-20 & 0.30 & 100.0\% & 0.0\% & 1 & 0.05 & 100.0\% & 0.0\% & 1 & 0.03 & -- & 100.0\% & 0.0\% & 1 & 0.03 & -- \\
a2-24 & 3.73 & 71.2\% & 0.0\% & 1 & 0.07 & 61.6\% & 0.0\% & 1 & 0.07 & -- & 60.1\% & 0.0\% & 1 & 0.07 & -- \\
a3-18 & 1.15 & 100.0\% & 0.0\% & 1 & 0.04 & 88.9\% & 0.0\% & 1 & 0.03 & -- & 88.9\% & 0.0\% & 1 & 0.04 & -- \\
a3-24 & 0.55 & 41.3\% & 0.0\% & 1 & 0.07 & 41.3\% & 0.0\% & 1 & 0.07 & -- & 41.3\% & 0.0\% & 1 & 0.07 & -- \\
a3-30 & 5.84 & 92.0\% & 0.0\% & 1 & 0.06 & 35.4\% & 0.0\% & 1 & 0.06 & -- & 34.6\% & 0.0\% & 1 & 0.06 & -- \\
a3-36 & 5.09 & 100.0\% & 0.0\% & 1 & 0.08 & 98.1\% & 0.0\% & 1 & 0.07 & -- & 98.1\% & 0.0\% & 1 & 0.08 & -- \\
a4-16 & 0.82 & 38.9\% & 0.0\% & 1 & 0.06 & 38.9\% & 0.0\% & 1 & 0.06 & -- & 38.9\% & 0.0\% & 1 & 0.06 & -- \\
a4-24 & 0.37 & 0.0\% & 0.0\% & 1 & 0.04 & 0.0\% & 0.0\% & 1 & 0.05 & -- & 0.0\% & 0.0\% & 1 & 0.05 & -- \\
a4-32 & 2.12 & 100.0\% & 0.0\% & 1 & 0.27 & 89.1\% & 0.0\% & 1 & 0.13 & -- & 89.1\% & 0.0\% & 1 & 0.16 & -- \\
a4-40 & 2.20 & 77.0\% & 0.0\% & 1 & 0.31 & 77.0\% & 0.0\% & 1 & 0.26 & -- & 77.0\% & 0.0\% & 1 & 0.27 & -- \\
a4-48 & 0.88 & 100.0\% & 0.0\% & 1 & 0.74 & 100.0\% & 0.0\% & 1 & 0.39 & -47\% & 100.0\% & 0.0\% & 1 & 0.41 & -45\% \\
a5-40 & 0.46 & 0.0\% & 0.0\% & 1 & 0.16 & 0.0\% & 0.0\% & 1 & 0.18 & -- & 0.0\% & 0.0\% & 1 & 0.16 & -- \\
a5-50 & 0.00 & -- & -- & 1 & 0.24 & -- & -- & 1 & 0.22 & -- & -- & -- & 1 & 0.25 & -- \\
a5-60 & 0.55 & 100.0\% & 0.0\% & 1 & 0.31 & 78.8\% & 0.0\% & 1 & 0.28 & -- & 78.8\% & 0.0\% & 1 & 0.30 & -- \\
a6-48 & 0.00 & -- & -- & 1 & 0.27 & -- & -- & 1 & 0.34 & -- & -- & -- & 1 & 0.32 & -- \\
a6-60 & 0.56 & 45.9\% & 0.0\% & 1 & 0.61 & 45.9\% & 0.0\% & 1 & 0.55 & -10\% & 45.9\% & 0.0\% & 1 & 0.56 & -8\% \\
a6-72 & 0.09 & 100.0\% & 0.0\% & 1 & 0.51 & 100.0\% & 0.0\% & 1 & 0.47 & -8\% & 100.0\% & 0.0\% & 1 & 0.47 & -8\% \\
a7-56 & 0.41 & 0.0\% & 0.0\% & 1 & 0.45 & 0.0\% & 0.0\% & 1 & 0.43 & -- & 0.0\% & 0.0\% & 1 & 0.35 & -- \\
a7-70 & 1.02 & 33.9\% & 0.0\% & 1 & 0.85 & 33.9\% & 0.0\% & 1 & 0.87 & +2\% & 33.9\% & 0.0\% & 1 & 0.96 & +13\% \\
a7-84 & 0.83 & 0.0\% & 0.0\% & 1 & 0.85 & 0.0\% & 0.0\% & 1 & 0.79 & -7\% & 0.0\% & 0.0\% & 1 & 0.97 & +14\% \\
a8-64 & 1.09 & 0.0\% & 0.0\% & 1 & 0.54 & 0.0\% & 0.0\% & 1 & 0.62 & +15\% & 0.0\% & 0.0\% & 1 & 0.59 & +9\% \\
a8-80 & 10.34 & 0.0\% & 0.0\% & 1 & 0.70 & 0.0\% & 0.0\% & 1 & 0.73 & +4\% & 0.0\% & 0.0\% & 1 & 0.76 & +9\% \\
a8-96 & 0.34 & 13.8\% & 0.0\% & 1 & 1.33 & 13.8\% & 0.0\% & 1 & 1.47 & +11\% & 13.8\% & 0.0\% & 1 & 1.49 & +12\% \\
\midrule
Total &  & 54.2\% & 0.0\% & 24 & 8.7 & 49.0\% & 0.0\% & 24 & 8.2 & -6\% & 48.9\% & 0.0\% & 24 & 8.5 & -2\% \\
\bottomrule
\end{tabular}
\par\vspace{3pt}
\begin{minipage}{\textwidth}\scriptsize
Notes: means over three seeds of LPgap, rGap, Nodes and CPU. Obj.: optimal value, or the best incumbent over all configurations and seeds if no configuration proved optimality. LPgap: gap of the LP bound to the objective value of the run. rGap: gap of the bound when Gurobi leaves the root node to branch, 0 if the run is solved at the root node; --: objective zero or root node not processed. Nodes: branch-and-bound nodes (1: solved at the root node, 0: solved in presolve without a root node). CPU: pure Gurobi time in seconds (graph and model building excluded); TL (x\%) = time limit of 7200~s reached in every seed, largest final gap. $\Delta$CPU: percentage change of CPU against base; --: both below 0.5~s or both at the limit. Total: sums of CPU and Nodes, means of the gaps, percentage change of the summed CPU.
\end{minipage}
\end{table}

\begin{table}[htbp]
\centering\scriptsize\setlength{\tabcolsep}{1.2pt}\renewcommand{\arraystretch}{0.9}
\caption{Base versus base + cuts on the b instances (relative detour~\eqref{eq:detour}).}
\label{tab:ridecut_LAEB_b_study_detour}
\begin{tabular}{l r rrrr rrrrr rrrrr}
\toprule
 & & \multicolumn{4}{c}{base} & \multicolumn{5}{c}{base + \eqref{eq:ridecut}} & \multicolumn{5}{c}{base + \eqref{eq:lridecut}} \\
\cmidrule(lr){3-6}\cmidrule(lr){7-11}\cmidrule(lr){12-16}
Inst. & Obj. & LPgap & rGap & Nodes & CPU & LPgap & rGap & Nodes & CPU & $\Delta$CPU & LPgap & rGap & Nodes & CPU & $\Delta$CPU \\
\midrule
b2-16 & 2.03 & 100.0\% & 0.0\% & 1 & 0.04 & 100.0\% & 0.0\% & 1 & 0.05 & -- & 100.0\% & 0.0\% & 1 & 0.05 & -- \\
b2-20 & 19.49 & 93.3\% & -- & 0 & 0.02 & 75.0\% & -- & 0 & 0.01 & -- & 75.0\% & -- & 0 & 0.01 & -- \\
b2-24 & 2.42 & 100.0\% & 0.0\% & 1 & 0.05 & 100.0\% & 0.0\% & 1 & 0.04 & -- & 100.0\% & 0.0\% & 1 & 0.04 & -- \\
b3-18 & 3.94 & 100.0\% & 0.0\% & 1 & 0.05 & 100.0\% & 0.0\% & 1 & 0.05 & -- & 100.0\% & 0.0\% & 1 & 0.06 & -- \\
b3-24 & 14.36 & 85.4\% & 0.0\% & 1 & 0.09 & 80.4\% & 0.0\% & 1 & 0.11 & -- & 80.4\% & 0.0\% & 1 & 0.12 & -- \\
b3-30 & 9.36 & 57.3\% & 0.0\% & 1 & 0.17 & 57.3\% & 0.0\% & 1 & 0.07 & -- & 57.3\% & 0.0\% & 1 & 0.08 & -- \\
b3-36 & 21.02 & 23.8\% & 0.0\% & 1 & 0.08 & 17.4\% & 0.0\% & 1 & 0.07 & -- & 17.4\% & 0.0\% & 1 & 0.08 & -- \\
b4-16 & 3.22 & 46.1\% & 0.0\% & 1 & 0.02 & 46.1\% & 0.0\% & 1 & 0.02 & -- & 46.1\% & 0.0\% & 1 & 0.02 & -- \\
b4-24 & 18.03 & 100.0\% & 0.0\% & 1 & 0.06 & 100.0\% & 0.0\% & 1 & 0.07 & -- & 100.0\% & 0.0\% & 1 & 0.08 & -- \\
b4-32 & 0.31 & 100.0\% & 0.0\% & 1 & 0.06 & 100.0\% & 0.0\% & 1 & 0.06 & -- & 100.0\% & 0.0\% & 1 & 0.06 & -- \\
b4-40 & 2.29 & 79.6\% & 0.0\% & 1 & 0.20 & 79.6\% & 0.0\% & 1 & 0.19 & -- & 79.6\% & 0.0\% & 1 & 0.16 & -- \\
b4-48 & 8.18 & 100.0\% & 0.0\% & 1 & 0.41 & 100.0\% & 0.0\% & 1 & 0.36 & -- & 100.0\% & 0.0\% & 1 & 0.35 & -- \\
b5-40 & 0.00 & -- & -- & 1 & 0.19 & -- & -- & 1 & 0.16 & -- & -- & -- & 1 & 0.17 & -- \\
b5-50 & 6.06 & 92.5\% & 0.0\% & 1 & 0.53 & 91.4\% & 0.0\% & 1 & 0.32 & -40\% & 91.4\% & 0.0\% & 1 & 0.34 & -36\% \\
b5-60 & 1.79 & 100.0\% & 0.0\% & 1 & 0.78 & 100.0\% & 0.0\% & 1 & 0.59 & -24\% & 100.0\% & 0.0\% & 1 & 0.69 & -12\% \\
b6-48 & 6.20 & 36.4\% & 0.0\% & 1 & 0.31 & 36.4\% & 0.0\% & 1 & 0.22 & -- & 36.4\% & 0.0\% & 1 & 0.38 & -- \\
b6-60 & 0.00 & -- & -- & 1 & 0.29 & -- & -- & 1 & 0.36 & -- & -- & -- & 1 & 0.33 & -- \\
b6-72 & 1.12 & 0.0\% & 0.0\% & 1 & 0.63 & 0.0\% & 0.0\% & 1 & 0.67 & +6\% & 0.0\% & 0.0\% & 1 & 0.74 & +17\% \\
b7-56 & 0.00 & -- & -- & 1 & 0.99 & -- & -- & 1 & 0.97 & -2\% & -- & -- & 1 & 0.95 & -4\% \\
b7-70 & 1.45 & 66.0\% & 0.0\% & 1 & 0.64 & 66.0\% & 0.0\% & 1 & 0.64 & +0\% & 66.0\% & 3.8\% & 2 & 0.73 & +14\% \\
b7-84 & 1.32 & 100.0\% & 14.3\% & 7 & 1.32 & 100.0\% & 0.0\% & 1 & 0.90 & -32\% & 100.0\% & 0.0\% & 1 & 1.07 & -19\% \\
b8-64 & 0.00 & -- & -- & 1 & 0.64 & -- & -- & 1 & 0.61 & -5\% & -- & -- & 1 & 0.53 & -17\% \\
b8-80 & 14.70 & 3.5\% & 0.0\% & 1 & 0.70 & 3.5\% & 0.0\% & 1 & 0.65 & -7\% & 3.5\% & 0.0\% & 1 & 0.65 & -7\% \\
b8-96 & 0.00 & -- & -- & 1 & 1.89 & -- & -- & 1 & 2.41 & +28\% & -- & -- & 1 & 2.32 & +23\% \\
\midrule
Total &  & 72.8\% & 0.8\% & 29 & 10.2 & 71.2\% & 0.0\% & 23 & 9.6 & -6\% & 71.2\% & 0.2\% & 24 & 10.0 & -2\% \\
\bottomrule
\end{tabular}
\par\vspace{3pt}
\begin{minipage}{\textwidth}\scriptsize
Notes: means over three seeds of LPgap, rGap, Nodes and CPU. Obj.: optimal value, or the best incumbent over all configurations and seeds if no configuration proved optimality. LPgap: gap of the LP bound to the objective value of the run. rGap: gap of the bound when Gurobi leaves the root node to branch, 0 if the run is solved at the root node; --: objective zero or root node not processed. Nodes: branch-and-bound nodes (1: solved at the root node, 0: solved in presolve without a root node). CPU: pure Gurobi time in seconds (graph and model building excluded); TL (x\%) = time limit of 7200~s reached in every seed, largest final gap. $\Delta$CPU: percentage change of CPU against base; --: both below 0.5~s or both at the limit. Total: sums of CPU and Nodes, means of the gaps, percentage change of the summed CPU.
\end{minipage}
\end{table}

\begin{table}[htbp]
\centering\scriptsize\setlength{\tabcolsep}{1.2pt}\renewcommand{\arraystretch}{0.9}
\caption{Base versus base + cuts on the a-X$^\star$ instances (relative detour~\eqref{eq:detour}).}
\label{tab:ridecut_LAEB_aX_study_detour}
\begin{tabular}{l r rrrr rrrrr rrrrr}
\toprule
 & & \multicolumn{4}{c}{base} & \multicolumn{5}{c}{base + \eqref{eq:ridecut}} & \multicolumn{5}{c}{base + \eqref{eq:lridecut}} \\
\cmidrule(lr){3-6}\cmidrule(lr){7-11}\cmidrule(lr){12-16}
Inst. & Obj. & LPgap & rGap & Nodes & CPU & LPgap & rGap & Nodes & CPU & $\Delta$CPU & LPgap & rGap & Nodes & CPU & $\Delta$CPU \\
\midrule
a2-16 & 1.27 & 100.0\% & 0.0\% & 1 & 0.03 & 100.0\% & 0.0\% & 1 & 0.07 & -- & 100.0\% & 0.0\% & 1 & 0.03 & -- \\
a2-20 & 0.30 & 100.0\% & 0.0\% & 1 & 0.06 & 100.0\% & 0.0\% & 1 & 0.05 & -- & 100.0\% & 0.0\% & 1 & 0.05 & -- \\
a2-24 & 3.73 & 100.0\% & 19.5\% & 19 & 0.48 & 100.0\% & 3.7\% & 6 & 0.23 & -- & 100.0\% & 3.7\% & 6 & 0.24 & -- \\
a3-18 & 1.15 & 100.0\% & 0.0\% & 1 & 0.12 & 100.0\% & 0.0\% & 1 & 0.09 & -- & 100.0\% & 0.0\% & 1 & 0.10 & -- \\
a3-24 & 0.55 & 41.3\% & 0.0\% & 1 & 0.23 & 41.3\% & 0.0\% & 1 & 0.21 & -- & 41.3\% & 0.0\% & 1 & 0.22 & -- \\
a3-30 & 4.19 & 100.0\% & 0.0\% & 1 & 0.33 & 100.0\% & 0.0\% & 1 & 0.14 & -- & 100.0\% & 0.0\% & 1 & 0.14 & -- \\
a3-36 & 5.09 & 100.0\% & 49.3\% & 520 & 1.30 & 100.0\% & 20.3\% & 43 & 0.85 & -35\% & 100.0\% & 16.0\% & 29 & 1.16 & -11\% \\
a4-16 & 0.82 & 38.9\% & 0.0\% & 1 & 0.20 & 38.9\% & 0.0\% & 1 & 0.18 & -- & 38.9\% & 0.0\% & 1 & 0.16 & -- \\
a4-24 & 0.37 & 0.0\% & 0.0\% & 1 & 0.15 & 0.0\% & 0.0\% & 1 & 0.13 & -- & 0.0\% & 0.0\% & 1 & 0.13 & -- \\
a4-32 & 2.12 & 100.0\% & 35.8\% & 119 & 1.29 & 100.0\% & 0.0\% & 1 & 0.53 & -59\% & 100.0\% & 0.0\% & 1 & 0.50 & -61\% \\
a4-40 & 2.20 & 77.0\% & 0.0\% & 1 & 1.35 & 77.0\% & 0.0\% & 1 & 0.97 & -28\% & 77.0\% & 0.0\% & 1 & 1.18 & -13\% \\
a4-48 & 0.88 & 100.0\% & 22.4\% & 362 & 2.87 & 100.0\% & 0.0\% & 1 & 2.03 & -29\% & 100.0\% & 0.0\% & 1 & 1.57 & -45\% \\
a5-40 & 0.46 & 0.0\% & 0.0\% & 1 & 0.49 & 0.0\% & 0.0\% & 1 & 0.48 & -- & 0.0\% & 0.0\% & 1 & 0.48 & -- \\
a5-50 & 0.00 & -- & -- & 1 & 0.95 & -- & -- & 1 & 0.72 & -24\% & -- & -- & 1 & 0.82 & -14\% \\
a5-60 & 0.55 & 100.0\% & 0.0\% & 1 & 1.26 & 100.0\% & 0.0\% & 1 & 1.01 & -20\% & 100.0\% & 0.0\% & 1 & 1.18 & -6\% \\
a6-48 & 0.00 & -- & -- & 1 & 1.71 & -- & -- & 1 & 1.35 & -21\% & -- & -- & 1 & 1.57 & -8\% \\
a6-60 & 0.56 & 45.9\% & 45.9\% & 260 & 6.93 & 45.9\% & 0.0\% & 1 & 4.83 & -30\% & 45.9\% & 0.0\% & 1 & 4.44 & -36\% \\
a6-72 & 0.09 & 100.0\% & 0.0\% & 1 & 2.25 & 100.0\% & 0.0\% & 1 & 1.96 & -13\% & 100.0\% & 0.0\% & 1 & 2.46 & +9\% \\
a7-56 & 0.41 & 0.0\% & 0.0\% & 1 & 2.01 & 0.0\% & 0.0\% & 1 & 1.63 & -19\% & 0.0\% & 0.0\% & 1 & 1.70 & -15\% \\
a7-70 & 1.02 & 33.9\% & 0.0\% & 1 & 3.42 & 33.9\% & 0.0\% & 1 & 3.22 & -6\% & 33.9\% & 0.0\% & 1 & 4.33 & +27\% \\
a7-84 & 0.83 & 0.0\% & 0.0\% & 1 & 5.31 & 0.0\% & 0.0\% & 1 & 3.96 & -25\% & 0.0\% & 0.0\% & 1 & 5.33 & +0\% \\
a8-64 & 1.09 & 0.0\% & 0.0\% & 1 & 4.97 & 0.0\% & 0.0\% & 1 & 3.50 & -30\% & 0.0\% & 0.0\% & 1 & 5.81 & +17\% \\
a8-80 & 10.34 & 0.0\% & 0.0\% & 1 & 5.18 & 0.0\% & 0.0\% & 1 & 4.67 & -10\% & 0.0\% & 0.0\% & 1 & 6.42 & +24\% \\
a8-96 & 0.34 & 13.8\% & 0.0\% & 1 & 10.98 & 13.8\% & 0.0\% & 1 & 9.96 & -9\% & 13.8\% & 0.0\% & 1 & 13.04 & +19\% \\
\midrule
Total &  & 56.9\% & 7.9\% & 1299 & 53.9 & 56.9\% & 1.1\% & 71 & 42.8 & -21\% & 56.9\% & 0.9\% & 57 & 53.0 & -2\% \\
\bottomrule
\end{tabular}
\par\vspace{3pt}
\begin{minipage}{\textwidth}\scriptsize
Notes: means over three seeds of LPgap, rGap, Nodes and CPU. Obj.: optimal value, or the best incumbent over all configurations and seeds if no configuration proved optimality. LPgap: gap of the LP bound to the objective value of the run. rGap: gap of the bound when Gurobi leaves the root node to branch, 0 if the run is solved at the root node; --: objective zero or root node not processed. Nodes: branch-and-bound nodes (1: solved at the root node, 0: solved in presolve without a root node). CPU: pure Gurobi time in seconds (graph and model building excluded); TL (x\%) = time limit of 7200~s reached in every seed, largest final gap. $\Delta$CPU: percentage change of CPU against base; --: both below 0.5~s or both at the limit. Total: sums of CPU and Nodes, means of the gaps, percentage change of the summed CPU.
\end{minipage}
\end{table}

\begin{table}[htbp]
\centering\scriptsize\setlength{\tabcolsep}{1.2pt}\renewcommand{\arraystretch}{0.9}
\caption{Base versus base + cuts on the b-X$^\star$ instances (relative detour~\eqref{eq:detour}).}
\label{tab:ridecut_LAEB_bX_study_detour}
\begin{tabular}{l r rrrr rrrrr rrrrr}
\toprule
 & & \multicolumn{4}{c}{base} & \multicolumn{5}{c}{base + \eqref{eq:ridecut}} & \multicolumn{5}{c}{base + \eqref{eq:lridecut}} \\
\cmidrule(lr){3-6}\cmidrule(lr){7-11}\cmidrule(lr){12-16}
Inst. & Obj. & LPgap & rGap & Nodes & CPU & LPgap & rGap & Nodes & CPU & $\Delta$CPU & LPgap & rGap & Nodes & CPU & $\Delta$CPU \\
\midrule
b2-16 & 2.03 & 100.0\% & 0.0\% & 1 & 0.11 & 100.0\% & 0.0\% & 1 & 0.09 & -- & 100.0\% & 0.0\% & 1 & 0.09 & -- \\
b2-20 & 18.50 & 100.0\% & 11.7\% & 11 & 0.14 & 100.0\% & 0.0\% & 1 & 0.04 & -- & 100.0\% & 0.0\% & 1 & 0.04 & -- \\
b2-24 & 2.42 & 100.0\% & 0.0\% & 1 & 0.10 & 100.0\% & 0.0\% & 1 & 0.09 & -- & 100.0\% & 0.0\% & 1 & 0.09 & -- \\
b3-18 & 2.43 & 100.0\% & 0.0\% & 1 & 0.11 & 100.0\% & 0.0\% & 1 & 0.10 & -- & 100.0\% & 0.0\% & 1 & 0.10 & -- \\
b3-24 & 10.37 & 79.9\% & 28.7\% & 103 & 1.00 & 79.9\% & 0.0\% & 1 & 0.61 & -39\% & 79.9\% & 0.0\% & 1 & 0.70 & -30\% \\
b3-30 & 7.24 & 87.9\% & 0.0\% & 1 & 0.45 & 87.9\% & 0.0\% & 1 & 0.23 & -- & 87.9\% & 0.0\% & 1 & 0.23 & -- \\
b3-36 & 12.70 & 100.0\% & 0.0\% & 1 & 0.22 & 100.0\% & 0.0\% & 1 & 0.21 & -- & 100.0\% & 0.0\% & 1 & 0.20 & -- \\
b4-16 & 3.22 & 46.1\% & 0.0\% & 1 & 0.04 & 46.1\% & 0.0\% & 1 & 0.04 & -- & 46.1\% & 0.0\% & 1 & 0.04 & -- \\
b4-24 & 12.08 & 100.0\% & 0.0\% & 1 & 0.17 & 100.0\% & 0.0\% & 1 & 0.16 & -- & 100.0\% & 0.0\% & 1 & 0.20 & -- \\
b4-32 & 0.31 & 100.0\% & 0.0\% & 1 & 0.12 & 100.0\% & 0.0\% & 1 & 0.11 & -- & 100.0\% & 0.0\% & 1 & 0.13 & -- \\
b4-40 & 2.29 & 79.6\% & 0.0\% & 1 & 0.45 & 79.6\% & 0.0\% & 1 & 0.31 & -- & 79.6\% & 0.0\% & 1 & 0.29 & -- \\
b4-48 & 7.31 & 100.0\% & 0.0\% & 1 & 2.02 & 100.0\% & 0.0\% & 1 & 1.13 & -44\% & 100.0\% & 0.0\% & 1 & 1.52 & -25\% \\
b5-40 & 0.00 & -- & -- & 1 & 0.65 & -- & -- & 1 & 0.51 & -22\% & -- & -- & 1 & 0.59 & -9\% \\
b5-50 & 6.06 & 92.5\% & 4.7\% & 14 & 2.06 & 92.5\% & 0.0\% & 1 & 1.21 & -41\% & 92.5\% & 2.1\% & 8 & 1.35 & -34\% \\
b5-60 & 1.79 & 100.0\% & 16.5\% & 22 & 2.43 & 100.0\% & 0.0\% & 1 & 1.88 & -23\% & 100.0\% & 0.0\% & 1 & 1.79 & -26\% \\
b6-48 & 6.20 & 100.0\% & 0.0\% & 1 & 0.59 & 100.0\% & 0.0\% & 1 & 0.48 & -19\% & 100.0\% & 0.0\% & 1 & 0.55 & -7\% \\
b6-60 & 0.00 & -- & -- & 1 & 0.68 & -- & -- & 1 & 0.79 & +16\% & -- & -- & 1 & 0.85 & +25\% \\
b6-72 & 1.12 & 0.0\% & 0.0\% & 1 & 2.51 & 0.0\% & 0.0\% & 1 & 1.89 & -25\% & 0.0\% & 0.0\% & 1 & 2.00 & -20\% \\
b7-56 & 0.00 & -- & -- & 1 & 6.68 & -- & -- & 1 & 5.18 & -22\% & -- & -- & 1 & 8.73 & +31\% \\
b7-70 & 1.45 & 66.0\% & 0.0\% & 1 & 2.87 & 66.0\% & 0.0\% & 1 & 2.18 & -24\% & 66.0\% & 0.0\% & 1 & 1.89 & -34\% \\
b7-84 & 1.32 & 100.0\% & 83.0\% & 301 & 4.54 & 100.0\% & 16.0\% & 14 & 3.65 & -20\% & 100.0\% & 7.7\% & 13 & 3.57 & -21\% \\
b8-64 & 0.00 & -- & -- & 1 & 2.90 & -- & -- & 1 & 1.79 & -38\% & -- & -- & 1 & 2.77 & -4\% \\
b8-80 & 14.70 & 3.5\% & 0.0\% & 1 & 2.36 & 3.5\% & 0.0\% & 1 & 2.25 & -5\% & 3.5\% & 0.0\% & 1 & 2.08 & -12\% \\
b8-96 & 0.00 & -- & -- & 1 & 15.71 & -- & -- & 1 & 12.36 & -21\% & -- & -- & 1 & 19.67 & +25\% \\
\midrule
Total &  & 81.9\% & 7.6\% & 470 & 48.9 & 81.9\% & 0.8\% & 37 & 37.3 & -24\% & 81.9\% & 0.5\% & 43 & 49.5 & +1\% \\
\bottomrule
\end{tabular}
\par\vspace{3pt}
\begin{minipage}{\textwidth}\scriptsize
Notes: means over three seeds of LPgap, rGap, Nodes and CPU. Obj.: optimal value, or the best incumbent over all configurations and seeds if no configuration proved optimality. LPgap: gap of the LP bound to the objective value of the run. rGap: gap of the bound when Gurobi leaves the root node to branch, 0 if the run is solved at the root node; --: objective zero or root node not processed. Nodes: branch-and-bound nodes (1: solved at the root node, 0: solved in presolve without a root node). CPU: pure Gurobi time in seconds (graph and model building excluded); TL (x\%) = time limit of 7200~s reached in every seed, largest final gap. $\Delta$CPU: percentage change of CPU against base; --: both below 0.5~s or both at the limit. Total: sums of CPU and Nodes, means of the gaps, percentage change of the summed CPU.
\end{minipage}
\end{table}

With this objective the benchmark instances are easy for every configuration: all 1440 solves end with a proven optimum, the longest in 25 seconds, and the optimal values are small, between 0 and 10.3 on the a sets and between 0 and 21.0 on the b sets, with an optimum of zero on two a and five b instances, where every request can be delivered at its earliest time.
The relative gaps in the tables are therefore large in percentage terms although small in absolute terms: the LP relaxation bound is zero on many instances, because the Big $M$ time constraints are not binding in the relaxation. The inequalities raise it only on some of the a and b instances without extended time windows, as they bound the delivery times only through the arcs on which a request is on board.

The effect of the inequalities \eqref{eq:ridecut} shows again in the search trees and in the root gaps.
On the a-X$^\star$ instances the base configuration branches on five instances and needs 1299 nodes and 53.9 seconds in total (mean root gap 7.9~percent), the configuration with the inequalities branches on two instances and needs 71 nodes and 42.8 seconds (mean root gap 1.1~percent); it is faster on all fifteen instances above half a second, by 24~percent on average.
On the b-X$^\star$ instances the totals fall from 470 nodes and 48.9 seconds to 37 nodes and 37.3 seconds and the mean root gap from 7.6 to 0.8~percent, faster on thirteen of fourteen instances, by 23~percent on average. The a and b instances are solved at the root node by all three configurations (b2-20 already in presolve, without any node), except b7-84 with the base configuration (7 nodes) and b7-70 with base + \eqref{eq:lridecut} (2 nodes), in about 9 and 10 seconds in total, and the inequalities~\eqref{eq:ridecut} change the totals by $-5$ and $-6$~percent.
The weaker configurations, run once each, gain more: on a-X$^\star$ the plain model needs 12\,392 nodes and 54.9 seconds and the plain model with the inequalities 193 nodes and 20.3 seconds, GP needs 1561 nodes and 31.7 seconds and GP with the inequalities 100 nodes and 19.1 seconds; on b-X$^\star$ the corresponding totals are 4352 nodes and 44.0 seconds against 435 nodes and 20.0 seconds for the plain model and 1651 nodes and 27.1 seconds against 61 nodes and 17.2 seconds for GP.

With this objective the strengthened inequalities~\eqref{eq:lridecut} are slower than~\eqref{eq:ridecut}: the totals are 53.0 against 42.8 seconds on a-X$^\star$ and 49.5 against 37.3 seconds on b-X$^\star$, which is the level of the base configuration, and 8.5 against 8.2 and 10.0 against 9.6 seconds on a and b; the weaker configurations lose one to five seconds each.
The search trees are as small as with~\eqref{eq:ridecut} (57 against 71 nodes on a-X$^\star$, 43 against 37 on b-X$^\star$) and the LP bounds are identical; the additional time is spent in the root node, where Gurobi's cutting-plane and heuristic loop runs longer with the coefficients $\tau_a$, while the LP relaxations take about the same time.

On the sets with windows of 45 minutes (Tables~\ref{tab:ridecut_LAEB_aX2_study_x2_detour} and~\ref{tab:ridecut_LAEB_bX2_study_x2_detour}) the relative-detour objective stays easy: all 432 solves end with a proven optimum, the longest in 156 seconds (a8-96-X$^{\star\star}$, base), with optimal values between 0 and 10.3 on a-X$^{\star\star}$ and between 0 and 18.5 on b-X$^{\star\star}$.

\begin{table}[htbp]
\centering\scriptsize\setlength{\tabcolsep}{1.2pt}\renewcommand{\arraystretch}{0.9}
\caption{Base versus base + cuts on the a-X$^{\star\star}$ instances (relative detour~\eqref{eq:detour}).}
\label{tab:ridecut_LAEB_aX2_study_x2_detour}
\begin{tabular}{l r rrrr rrrrr rrrrr}
\toprule
 & & \multicolumn{4}{c}{base} & \multicolumn{5}{c}{base + \eqref{eq:ridecut}} & \multicolumn{5}{c}{base + \eqref{eq:lridecut}} \\
\cmidrule(lr){3-6}\cmidrule(lr){7-11}\cmidrule(lr){12-16}
Inst. & Obj. & LPgap & rGap & Nodes & CPU & LPgap & rGap & Nodes & CPU & $\Delta$CPU & LPgap & rGap & Nodes & CPU & $\Delta$CPU \\
\midrule
a2-16 & 1.27 & 100.0\% & 0.0\% & 1 & 0.06 & 100.0\% & 0.0\% & 1 & 0.04 & -- & 100.0\% & 0.0\% & 1 & 0.05 & -- \\
a2-20 & 0.30 & 100.0\% & 0.0\% & 1 & 0.14 & 100.0\% & 0.0\% & 1 & 0.12 & -- & 100.0\% & 0.0\% & 1 & 0.10 & -- \\
a2-24 & 3.73 & 100.0\% & 38.6\% & 74 & 0.60 & 100.0\% & 0.0\% & 1 & 0.44 & -27\% & 100.0\% & 10.1\% & 15 & 0.50 & -17\% \\
a3-18 & 1.15 & 100.0\% & 0.0\% & 1 & 0.29 & 100.0\% & 0.0\% & 1 & 0.32 & -- & 100.0\% & 0.0\% & 1 & 0.21 & -- \\
a3-24 & 0.55 & 41.3\% & 0.0\% & 1 & 0.42 & 41.3\% & 0.0\% & 1 & 0.39 & -- & 41.3\% & 0.0\% & 1 & 0.35 & -- \\
a3-30 & 4.19 & 100.0\% & 49.1\% & 112 & 1.23 & 100.0\% & 0.0\% & 1 & 0.30 & -76\% & 100.0\% & 0.0\% & 1 & 0.31 & -75\% \\
a3-36 & 4.87 & 100.0\% & 72.3\% & 1545 & 5.40 & 100.0\% & 19.4\% & 91 & 1.97 & -64\% & 100.0\% & 25.0\% & 115 & 1.79 & -67\% \\
a4-16 & 0.82 & 38.9\% & 0.0\% & 1 & 0.43 & 38.9\% & 0.0\% & 1 & 0.37 & -- & 38.9\% & 0.0\% & 1 & 0.35 & -- \\
a4-24 & 0.37 & 0.0\% & 0.0\% & 1 & 0.46 & 0.0\% & 0.0\% & 1 & 0.31 & -- & 0.0\% & 0.0\% & 1 & 0.38 & -- \\
a4-32 & 2.12 & 100.0\% & 60.5\% & 515 & 3.30 & 100.0\% & 0.0\% & 1 & 1.17 & -65\% & 100.0\% & 0.0\% & 1 & 1.34 & -59\% \\
a4-40 & 2.20 & 77.0\% & 0.0\% & 1 & 5.00 & 77.0\% & 0.0\% & 1 & 3.05 & -39\% & 77.0\% & 0.0\% & 1 & 4.50 & -10\% \\
a4-48 & 0.88 & 100.0\% & 65.6\% & 266 & 10.91 & 100.0\% & 0.0\% & 1 & 6.22 & -43\% & 100.0\% & 0.0\% & 1 & 7.35 & -33\% \\
a5-40 & 0.46 & 0.0\% & 0.0\% & 1 & 1.30 & 0.0\% & 0.0\% & 1 & 1.23 & -5\% & 0.0\% & 0.0\% & 1 & 1.22 & -6\% \\
a5-50 & 0.00 & -- & -- & 1 & 2.80 & -- & -- & 1 & 2.20 & -21\% & -- & -- & 1 & 2.61 & -7\% \\
a5-60 & 0.55 & 100.0\% & 0.0\% & 1 & 5.33 & 100.0\% & 0.0\% & 1 & 4.08 & -23\% & 100.0\% & 0.0\% & 1 & 4.27 & -20\% \\
a6-48 & 0.00 & -- & -- & 1 & 12.57 & -- & -- & 1 & 11.26 & -10\% & -- & -- & 1 & 11.52 & -8\% \\
a6-60 & 0.56 & 45.9\% & 30.6\% & 301 & 37.05 & 45.9\% & 0.0\% & 1 & 28.15 & -24\% & 45.9\% & 0.0\% & 1 & 21.94 & -41\% \\
a6-72 & 0.09 & 100.0\% & 0.0\% & 1 & 18.57 & 100.0\% & 0.0\% & 1 & 11.76 & -37\% & 100.0\% & 0.0\% & 1 & 13.27 & -29\% \\
a7-56 & 0.41 & 0.0\% & 0.0\% & 1 & 11.48 & 0.0\% & 0.0\% & 1 & 10.91 & -5\% & 0.0\% & 0.0\% & 1 & 10.75 & -6\% \\
a7-70 & 1.02 & 33.9\% & 0.0\% & 1 & 34.74 & 33.9\% & 0.0\% & 1 & 24.25 & -30\% & 33.9\% & 0.0\% & 1 & 20.12 & -42\% \\
a7-84 & 0.83 & 0.0\% & 0.0\% & 3 & 50.60 & 0.0\% & 0.0\% & 1 & 25.46 & -50\% & 0.0\% & 0.0\% & 1 & 25.90 & -49\% \\
a8-64 & 1.09 & 0.0\% & 0.0\% & 1 & 44.40 & 0.0\% & 0.0\% & 1 & 25.37 & -43\% & 0.0\% & 0.0\% & 1 & 33.49 & -25\% \\
a8-80 & 10.34 & 0.0\% & 0.0\% & 1 & 49.09 & 0.0\% & 0.0\% & 1 & 34.10 & -31\% & 0.0\% & 0.0\% & 1 & 35.12 & -28\% \\
a8-96 & 0.34 & 13.8\% & 0.0\% & 1 & 144.81 & 13.8\% & 0.0\% & 1 & 69.95 & -52\% & 13.8\% & 0.0\% & 1 & 73.57 & -49\% \\
\midrule
Total &  & 56.9\% & 14.4\% & 2833 & 441.0 & 56.9\% & 0.9\% & 114 & 263.4 & -40\% & 56.9\% & 1.6\% & 152 & 271.0 & -39\% \\
\bottomrule
\end{tabular}
\par\vspace{3pt}
\begin{minipage}{\textwidth}\scriptsize
Notes: means over three seeds of LPgap, rGap, Nodes and CPU. Obj.: optimal value, or the best incumbent over all configurations and seeds if no configuration proved optimality. LPgap: gap of the LP bound to the objective value of the run. rGap: gap of the bound when Gurobi leaves the root node to branch, 0 if the run is solved at the root node; --: objective zero or root node not processed. Nodes: branch-and-bound nodes (1: solved at the root node, 0: solved in presolve without a root node). CPU: pure Gurobi time in seconds (graph and model building excluded); TL (x\%) = time limit of 7200~s reached in every seed, largest final gap. $\Delta$CPU: percentage change of CPU against base; --: both below 0.5~s or both at the limit. Total: sums of CPU and Nodes, means of the gaps, percentage change of the summed CPU.
\end{minipage}
\end{table}

\begin{table}[htbp]
\centering\scriptsize\setlength{\tabcolsep}{1.2pt}\renewcommand{\arraystretch}{0.9}
\caption{Base versus base + cuts on the b-X$^{\star\star}$ instances (relative detour~\eqref{eq:detour}).}
\label{tab:ridecut_LAEB_bX2_study_x2_detour}
\begin{tabular}{l r rrrr rrrrr rrrrr}
\toprule
 & & \multicolumn{4}{c}{base} & \multicolumn{5}{c}{base + \eqref{eq:ridecut}} & \multicolumn{5}{c}{base + \eqref{eq:lridecut}} \\
\cmidrule(lr){3-6}\cmidrule(lr){7-11}\cmidrule(lr){12-16}
Inst. & Obj. & LPgap & rGap & Nodes & CPU & LPgap & rGap & Nodes & CPU & $\Delta$CPU & LPgap & rGap & Nodes & CPU & $\Delta$CPU \\
\midrule
b2-16 & 2.03 & 100.0\% & 0.0\% & 1 & 0.33 & 100.0\% & 0.0\% & 1 & 0.23 & -- & 100.0\% & 0.0\% & 1 & 0.22 & -- \\
b2-20 & 18.50 & 100.0\% & 58.9\% & 101 & 0.29 & 100.0\% & 25.1\% & 16 & 0.25 & -- & 100.0\% & 25.9\% & 19 & 0.26 & -- \\
b2-24 & 2.42 & 100.0\% & 0.0\% & 1 & 0.16 & 100.0\% & 0.0\% & 1 & 0.18 & -- & 100.0\% & 0.0\% & 1 & 0.14 & -- \\
b3-18 & 2.43 & 100.0\% & 0.0\% & 1 & 0.23 & 100.0\% & 0.0\% & 1 & 0.18 & -- & 100.0\% & 0.0\% & 1 & 0.19 & -- \\
b3-24 & 10.37 & 79.9\% & 79.5\% & 3205 & 11.93 & 79.9\% & 53.0\% & 1226 & 3.02 & -75\% & 79.9\% & 53.4\% & 1134 & 2.79 & -77\% \\
b3-30 & 7.24 & 87.9\% & 0.0\% & 1 & 0.86 & 87.9\% & 0.0\% & 1 & 0.42 & -51\% & 87.9\% & 0.0\% & 1 & 0.44 & -49\% \\
b3-36 & 10.85 & 100.0\% & 0.0\% & 1 & 0.50 & 100.0\% & 0.0\% & 1 & 0.38 & -24\% & 100.0\% & 0.0\% & 1 & 0.37 & -26\% \\
b4-16 & 3.22 & 46.1\% & 0.0\% & 1 & 0.07 & 46.1\% & 0.0\% & 1 & 0.06 & -- & 46.1\% & 0.0\% & 1 & 0.06 & -- \\
b4-24 & 12.08 & 100.0\% & 0.0\% & 1 & 0.65 & 100.0\% & 0.0\% & 1 & 0.58 & -11\% & 100.0\% & 0.0\% & 1 & 0.69 & +6\% \\
b4-32 & 0.31 & 100.0\% & 0.0\% & 1 & 0.22 & 100.0\% & 0.0\% & 1 & 0.19 & -- & 100.0\% & 0.0\% & 1 & 0.17 & -- \\
b4-40 & 2.29 & 79.6\% & 0.0\% & 1 & 0.54 & 79.6\% & 0.0\% & 1 & 0.48 & -11\% & 79.6\% & 0.0\% & 1 & 0.57 & +6\% \\
b4-48 & 7.31 & 100.0\% & 0.0\% & 1 & 7.03 & 100.0\% & 0.0\% & 1 & 5.82 & -17\% & 100.0\% & 0.0\% & 1 & 8.14 & +16\% \\
b5-40 & 0.00 & -- & -- & 1 & 2.18 & -- & -- & 1 & 2.01 & -8\% & -- & -- & 1 & 1.79 & -18\% \\
b5-50 & 6.06 & 92.5\% & 69.8\% & 2040 & 8.62 & 92.5\% & 0.0\% & 1 & 2.74 & -68\% & 92.5\% & 25.3\% & 144 & 3.51 & -59\% \\
b5-60 & 1.79 & 100.0\% & 86.7\% & 102 & 7.98 & 100.0\% & 0.0\% & 1 & 5.00 & -37\% & 100.0\% & 0.0\% & 1 & 4.99 & -37\% \\
b6-48 & 6.20 & 100.0\% & 0.0\% & 1 & 2.48 & 100.0\% & 0.0\% & 1 & 1.29 & -48\% & 100.0\% & 0.0\% & 1 & 1.65 & -33\% \\
b6-60 & 0.00 & -- & -- & 1 & 1.66 & -- & -- & 1 & 1.98 & +19\% & -- & -- & 1 & 2.04 & +23\% \\
b6-72 & 1.12 & 0.0\% & 0.0\% & 1 & 10.34 & 0.0\% & 0.0\% & 1 & 7.15 & -31\% & 0.0\% & 0.0\% & 1 & 8.22 & -21\% \\
b7-56 & 0.00 & -- & -- & 1 & 61.23 & -- & -- & 1 & 10.73 & -82\% & -- & -- & 1 & 20.47 & -67\% \\
b7-70 & 1.45 & 66.0\% & 22.0\% & 29 & 8.54 & 66.0\% & 3.5\% & 2 & 4.36 & -49\% & 66.0\% & 0.0\% & 1 & 5.79 & -32\% \\
b7-84 & 1.32 & 100.0\% & 89.1\% & 375 & 11.07 & 100.0\% & 0.0\% & 1 & 11.54 & +4\% & 100.0\% & 3.4\% & 3 & 10.69 & -3\% \\
b8-64 & 0.00 & -- & -- & 1 & 13.64 & -- & -- & 1 & 8.46 & -38\% & -- & -- & 1 & 7.65 & -44\% \\
b8-80 & 14.70 & 3.5\% & 1.2\% & 3 & 6.17 & 3.5\% & 0.0\% & 1 & 4.89 & -21\% & 3.5\% & 0.0\% & 1 & 4.77 & -23\% \\
b8-96 & 0.00 & -- & -- & 1 & 108.48 & -- & -- & 1 & 42.03 & -61\% & -- & -- & 1 & 53.02 & -51\% \\
\midrule
Total &  & 81.9\% & 21.4\% & 5872 & 265.2 & 81.9\% & 4.3\% & 1265 & 114.0 & -57\% & 81.9\% & 5.7\% & 1320 & 138.6 & -48\% \\
\bottomrule
\end{tabular}
\par\vspace{3pt}
\begin{minipage}{\textwidth}\scriptsize
Notes: means over three seeds of LPgap, rGap, Nodes and CPU. Obj.: optimal value, or the best incumbent over all configurations and seeds if no configuration proved optimality. LPgap: gap of the LP bound to the objective value of the run. rGap: gap of the bound when Gurobi leaves the root node to branch, 0 if the run is solved at the root node; --: objective zero or root node not processed. Nodes: branch-and-bound nodes (1: solved at the root node, 0: solved in presolve without a root node). CPU: pure Gurobi time in seconds (graph and model building excluded); TL (x\%) = time limit of 7200~s reached in every seed, largest final gap. $\Delta$CPU: percentage change of CPU against base; --: both below 0.5~s or both at the limit. Total: sums of CPU and Nodes, means of the gaps, percentage change of the summed CPU.
\end{minipage}
\end{table}

On a-X$^{\star\star}$ the base configuration branches on seven instances and needs 2833 nodes and 441 seconds in total (mean root gap 14.4~percent), the configuration with the inequalities \eqref{eq:ridecut} branches on one instance and needs 114 nodes and 263 seconds (mean root gap 0.9~percent); it is faster on all eighteen instances above half a second, by 36~percent on average.
On b-X$^{\star\star}$ the totals fall from 5872 nodes and 265 seconds to 1265 nodes and 114 seconds and the mean root gap from 21.4 to 4.3~percent, faster on sixteen of eighteen instances, by 34~percent on average.
The strengthened inequalities~\eqref{eq:lridecut} need 271 and 139 seconds on the two sets, again slightly more than~\eqref{eq:ridecut}, with search trees of the same size.
The total gains by the new inequalities for this objective thus grow with the width of the time windows, from about 21-24\% on the X$^\star$ sets to about 40-57\% on the X$^{\star\star}$ sets. The strengthened versions do not have a positive effect here. One reason is probably that the objective function aims at avoiding waiting times at least with customers on board, which is not the case of routing cost is minimised.

\bibliographystyle{plainnat}
\bibliography{references}

\end{document}